%% file: main.tex
\documentclass[12pt]{amsart}
\usepackage{amsfonts}
\usepackage{amssymb}
\usepackage{graphicx}
\usepackage{fullpage}

\newtheorem{thm}{Theorem}[section]
\newtheorem{prop}[thm]{Proposition}
\newtheorem{cor}[thm]{Corollary}
\newtheorem{lemma}[thm]{Lemma}
\newtheorem{bigthm}{Theorem}

\newtheorem*{thmB*}{Theorem B$^*$}
\newtheorem*{thmB**}{Theorem B$^{**}$}

\theoremstyle{remark}
\newtheorem{ex}{Example}[section]

\def\C{\mathbb{C}}
\def\R{\mathbb{R}}
\def\N{\mathbb{N}}
\def\Z{\mathbb{Z}}
\def\Bc{\mathcal{B}}
\def\Ec{\mathcal{E}}
\def\Rc{\mathcal{R}}
\def\eps{\varepsilon}
\def\phi{\varphi}
\def\d{\partial}

\title{The external boundary of the bifurcation locus~$M_2$}
\author{V. Timorin}

\begin{document}

\begin{abstract}
Consider a quadratic rational self-map of the Riemann sphere such that
one critical point is periodic of period 2, and the other critical point
lies on the boundary of its immediate basin of attraction.
We will give explicit topological models for all such maps.
We also discuss the corresponding parameter picture.
\end{abstract}

\maketitle

\thispagestyle{empty}
\input{imsmark}
\SBIMSMark{2007/1}{February 2007}{}

\section{Introduction}
\label{s_intro}

\subsection{The family $V_2$}
Consider the set $V_2$ of holomorphic conjugacy classes of
quadratic rational maps that have a super-attracting periodic cycle of period 2
(we follow the notation of Mary Rees).
The complement in $V_2$ to the class of the single map $z\mapsto 1/z^2$
is denoted by $V_{2,0}$.
The set $V_{2,0}$ is parameterized by a single complex number.
Indeed, for any map $f$ in $V_{2,0}$, the critical point of period two can be
mapped to $\infty$, its $f$-image to $0$, and the other critical point to $-1$.
Then we obtain a map of the form
$$
f_a(z)=\frac{a}{z^2+2z},\quad a\ne 0
$$
holomorphically conjugate to $f$.
Thus the set $V_{2,0}$ is identified with $\C-0$.

The family $V_2$ is just the second term in the sequence $V_1, V_2, V_3,\dots$,
where, by definition, $V_n$ consists of holomorphic conjugacy classes
of quadratic rational maps with a periodic critical orbit of period $n$.
Such maps have one ``free'' critical point, hence each family $V_n$
has complex dimension 1.
Note that $V_1$ is the family of quadratic polynomials, i.e., holomorphic
endomorphisms of the Riemann sphere of degree 2 with a fixed critical point at $\infty$.
Any quadratic polynomial is holomorphically conjugate to a map $z\mapsto z^2+c$.
Thus $V_1$ can be identified with the complex $c$-plane.
For a map $z\mapsto z^2+c$, the ``free'' critical point is $0$.
The family $V_1$ is the most studied family in complex dynamics.
The main object describing the structure of $V_1$ is the {\em Mandelbrot set} $M$
defined as the set of all parameter values $c$ such that the orbit of the
critical point $0$ under $z\mapsto z^2+c$ is bounded.

Similarly to the case of quadratic polynomials, we can define the set $M_2$
(an analog of the Mandelbrot set for $V_2$) as the set of all parameter values $a$
such that the orbit of $-1$ under $f_a$ is bounded.
A conjectural description of the topology of $M_2$ is given in \cite{Wittner}.
In this paper, we deal with maps on the {\em external boundary} of $M_2$, i.e.
the boundary of the only unbounded component of $\C-M_2$.

\begin{figure}
\centering
\includegraphics{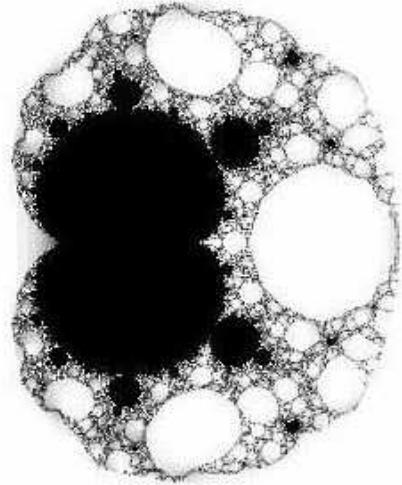}
\caption{The set $M_2$}
\end{figure}

In \cite{ReesV3}, M. Rees studies the parameter plane of $V_3$, which turns
out to be much more complicated than $V_2$.

\subsection{Invariant laminations}
Invariant laminations were introduced by Thurston \cite{Thurston} to describe
quadratic polynomials with locally connected Julia sets.
A set $L$ of hyperbolic geodesics in the open unit disk is a
{\em geodesic lamination} if any two different geodesics in $L$ do not intersect,
and the union of $L$ is closed with respect to the induced topology
on the unit disk.
For any pair of points $z$, $w$ on the unit circle, the geodesic with endpoints $z$ and $w$
will be written as $zw$.
Any geodesic lamination $L$ defines an equivalence relation $\sim_L$ on the unit circle $S^1$.
Namely, two different points on $S^1$ are equivalent if they are connected
by a leaf of $L$ or by a broken line consisting of leaves.
For many quadratic polynomials, the Julia set is homeomorphic to the
quotient of the unit circle by an equivalence relation $\sim_L$.

We say that a geodesic lamination $L$ on the unit circle
is {\em invariant under the map $z\mapsto z^2$} if the following conditions hold:
\begin{itemize}
\item if $z_1z_2\in L$, then $z_1^2z_2^2\in L$ or $z_1^2=z_2^2$,
\item if $z_1z_2\in L$, then $(-z_1)(-z_2)\in L$,
\item if $z_1^2z_2^2\in L$, then $z_1z_2\in L$ or $z_1(-z_2)\in L$.
\end{itemize}
Such laminations are also known as {\em quadratic invariant laminations}.
Any quadratic polynomial $p$ defines a quadratic invariant lamination.
In many cases, the quotient of the unit circle by the corresponding equivalence relation
is homeomorphic to the Julia set $J$, and the projection of $S^1$ onto $J$
semi-conjugates the map $z\mapsto z^2$ with the restriction of $p$ to $J$.

A {\em gap} of a geodesic lamination is any component of the complement
to all leaves in the unit disk.
Let $L$ be a quadratic invariant lamination.
The map $z\mapsto z^2$ admits a natural extension over all leaves and gaps of $L$.
This extension is called the {\em lamination map} of $L$ and is denoted by $s_L$.
The image of any leaf under $s_L$ is a leaf or a single point.
The image of any gap is a gap, or a leaf, or a single point.
Suppose that $L$ is {\em clean}, i.e. any two adjacent leaves of $L$
are sides of a common finite-sided gap.
Then we can also extend the equivalence relation $\sim_L$ to $\C$.
The equivalence classes of $\sim_L$ are defined as finite-sided gaps, leaves, or points.

If $L$ is clean, then the quotient $\C/\sim_L$ is homeomorphic to $\C$.
The lamination map $s_L$ defines a continuous self-map $[s_L]$ of this quotient.
We say that the lamination $L$ {\em models} a quadratic polynomial $p$ if
the quotient $\C/\sim_L$ is homeomorphic to $\C$, and the map $[s_L]$
is topologically conjugate to $p$.
E.g. any critically finite quadratic polynomial is modeled
by the corresponding quadratic invariant lamination.
The same is true for many quadratic polynomials with Siegel disks, but not
for quadratic polynomials with Cremer points.

Let $y_0$ be a real number between 0 and 1.
Denote by $l_0$ the diameter connecting the points $e^{\pi iy_0}$ and $-e^{\pi iy_0}$
on the unit circle.
Consider all geodesics $z_1z_2$ in the unit disk such that, for every $k$,
the geodesic $z_1^{2^k}z_2^{2^k}$ does not intersect $l_0$ or coincides with $l_0$.
This set of geodesics is an invariant lamination, which we denote by $L(y_0)$.
If a quadratic polynomial $p$ is modeled by $L(y_0)$, then $p$ belongs to the boundary of the Mandelbrot set.
There is a natural {\em parameter equivalence} relation on the unit circle.
Points $e^{2\pi iy_0}$ and $e^{2\pi iy'_0}$ are parameter equivalent if the laminations
$L(y_0)$ and $L(y'_0)$ correspond to the same quadratic polynomial in a certain
well-defined sense, although they may not model this polynomial (e.g. $L(0)$
corresponds to the parabolic map $z\mapsto z^2+1/4$, but the equivalence relation
$\sim_{L(0)}$ identifies all binary rational points on the unit circle).
It turns out that the parameter equivalence relation also corresponds to
a geodesic lamination in the unit disk.
This lamination is called the {\em parameter lamination}, or the {\em quadratic minor lamination}.
Thurston \cite{Thurston} gave a description of the parameter lamination using his ``minor leaf theory''.
Conjecturally, the boundary of the Mandelbrot set is homeomorphic to
the quotient of the unit circle by the parameter equivalence relation.
This conjecture is equivalent to the MLC conjecture (stating that the Mandelbrot
set is locally connected).

\subsection{Two-sided laminations}
In the theory of quadratic invariant laminations, the single quadratic polynomial
$z\mapsto z^2$ is used to build models for the dynamics of many other quadratic polynomials.
The Julia set of $z\mapsto z^2$ is the unit circle, and the unit disk is preserved.
A similar idea can be used to build models for rational maps of class $V_2$.
To this end, one can use the rational map $z\mapsto 1/z^2$.
This is the only map in $V_2$ not conjugate to a map of the form $f_a$.
Its Julia set is also the unit circle.
However, the map $z\mapsto 1/z^2$ interchanges the inside and the outside of
the unit disk.

Let us define an analog of quadratic invariant laminations for the map $z\mapsto 1/z^2$.
A {\em two-sided geodesic lamination} is a set of geodesics that live
both inside and outside of the unit disk.
Note that the outside of the unit disk is also a topological disk in $\overline{\C}$.
Geodesics are in the sense of the Poincar\'e metric (on the inside or on the outside of the unit disk).
We will sometimes use $2L$ to denote a two-sided lamination,
but this notation does not assume any multiplication by 2 (in other words,
$2L$ is to be thought of as a single piece of notation).
A two-sided lamination $2L$ gives rise to a pair of laminations
$L_0$ and $L_\infty$, where the leaves of $L_0$ are inside of the unit circle,
and the leaves of $L_\infty$ are outside.
The two-sided lamination $2L=L_0\cup L_\infty$ is called {\em invariant} under $z\mapsto 1/z^2$ if
the following conditions hold:
\begin{itemize}
\item if $z_1z_2\in L_0$, then $(1/z_1^2)(1/z_2^2)\in L_\infty$ or $z_1^2=z_2^2$,
\item if $z_1z_2\in L_0$, then $(-z_1)(-z_2)\in L_0$,
\item if $z_1^2z_2^2\in L_0$, then $z_1z_2\in L_\infty$ or $z_1(-z_2)\in L_\infty$,
\end{itemize}
and the same conditions with $L_0$ and $L_\infty$ interchanged.
Let $\sim_0$ and $\sim_\infty$ denote the equivalence relations on the
unit circle corresponding to the laminations $L_0$ and $L_\infty$, respectively.

Two-sided laminations were first considered by D. Ahmadi \cite{Ahmadi}.
He used a different language (``laminations on two disks'').
In \cite{Ahmadi}, a classification of two-sided laminations is given,
similar to the ``minor leaf theory'' of Thurston \cite{Thurston}.

Gaps of two-sided laminations and the corresponding lamination maps
are defined in the same way as for invariant laminations of the unit disk.
For a two-sided lamination $2L$, extend the equivalence relations $\sim_0$
and $\sim_\infty$ to the unit disk and to the outside of the unit disk,
respectively, in the same way as for invariant quadratic laminations.
Define $\sim_{2L}$ to be the smallest equivalence relation containing
both $\sim_0$ and $\sim_\infty$.
We say that $2L$ {\em models} a quadratic rational map $f$
if the quotient $\overline{\C}/\sim_{2L}$ is homeomorphic
to the sphere, and the map $[s_{2L}]$ is topologically conjugate to $f$.

We will now define a particular family of two-sided laminations invariant
under $z\mapsto 1/z^2$.
Let $x_0$ be a real number strictly between 0 and 1.
Consider the arc $\sigma_0$ of the unit circle bounded by the points $e^{2\pi ix_0}$ and
$-e^{2\pi ix_0}$ and not containing the point 1.
Let $\sigma$ be any component of the full $n$-fold preimage of $\sigma_0$ under $z\mapsto 1/z^2$.
Connect the endpoints of $\sigma$ by a geodesic in the complement to the unit circle.
This geodesic should be inside the unit circle if $n$ is even, and outside if $n$ is odd.
For certain values of $x_0$ (which we will describe explicitly later), the
set of geodesics thus constructed is a two-sided lamination.
We denote this lamination by $2L(x_0)$.
If $2L(x_0)$ exists, then it is clearly invariant under the map $z\mapsto 1/z^2$.

\subsection{Statement of the main theorems}
For a map $f_a\in V_2$, denote by $\Omega$ the immediate basin of attraction
of the critical cycle $\{0,\infty\}$.

\begin{bigthm}
\label{lc}
Suppose that $-1\in\d\Omega$.
Then the Julia set of $f_a$ is locally connected.
\end{bigthm}

Let $\Omega_0$ and $\Omega_\infty$ denote the components of $\Omega$
containing $0$ and $\infty$, respectively.
As we will see, the critical point $-1$ cannot be on the boundary of $\Omega_\infty$.
Thus, under the assumptions of Theorem \ref{lc}, we can only have $-1\in\d\Omega_0$.
We will prove in this case that $\overline{\Omega}_0$ is a closed topological disk.
Moreover, there is a homeomorphism $H$ of the closed unit disk to $\overline{\Omega}_0$
that conjugates the map $z\mapsto z^2$ with the map $f_a^{\circ 2}$.
We say that a point in $\overline{\Omega}_0$ {\em has angle $\theta$} if
this point coincides with $H(re^{2\pi i\theta})$ for some $0\le r\le 1$.

\begin{bigthm}
\label{bndry}
Suppose that the critical point $-1$ belongs to $\d\Omega_0$ and has angle $\theta_0$.
Then, for
$$
x_0=\sum_{m=1}^\infty\frac{[(2^m-1)\theta_0]+1}{2^{2m+1}},
$$
the two-sided lamination $2L(x_0)$ exists and models the map $f_a$.
\end{bigthm}

The maps $f_a$ from Theorems \ref{lc} and \ref{bndry}, together with countably many parabolic
maps, form the external boundary of $M_2$ (the boundary of the unbounded component of $\C-M_2$).
A more detailed statement will be given below.

\begin{figure}
\centering
\includegraphics[width=5cm]{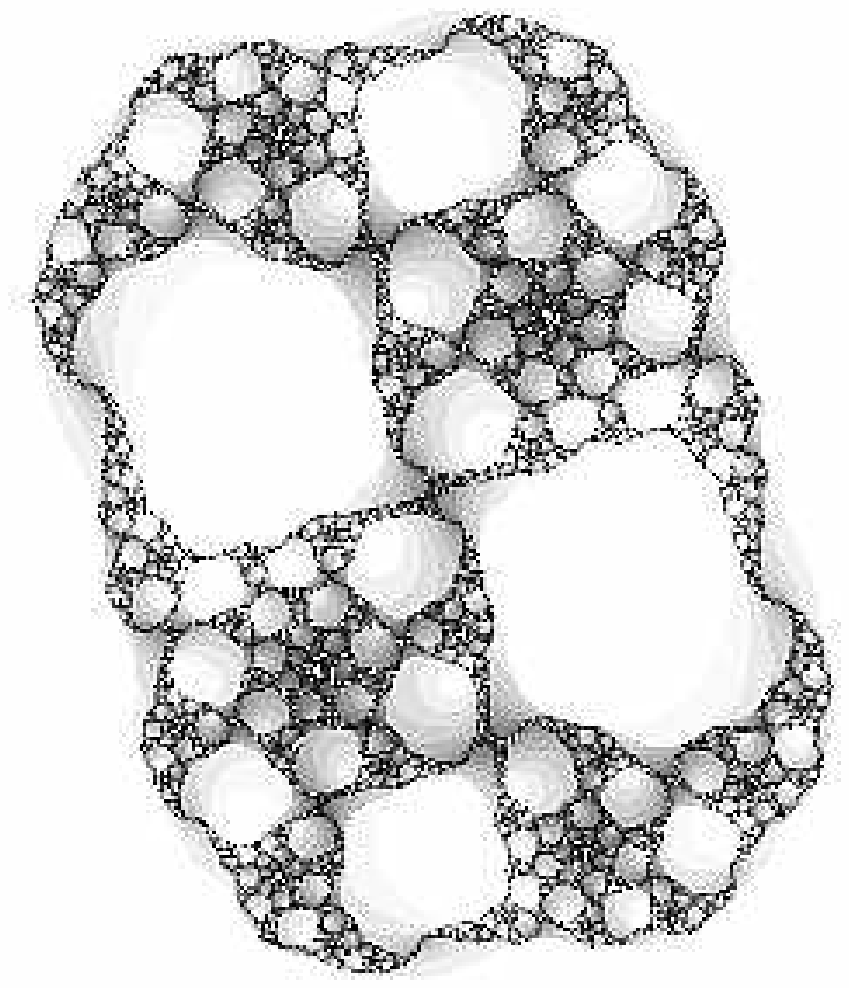}
\includegraphics[width=5cm]{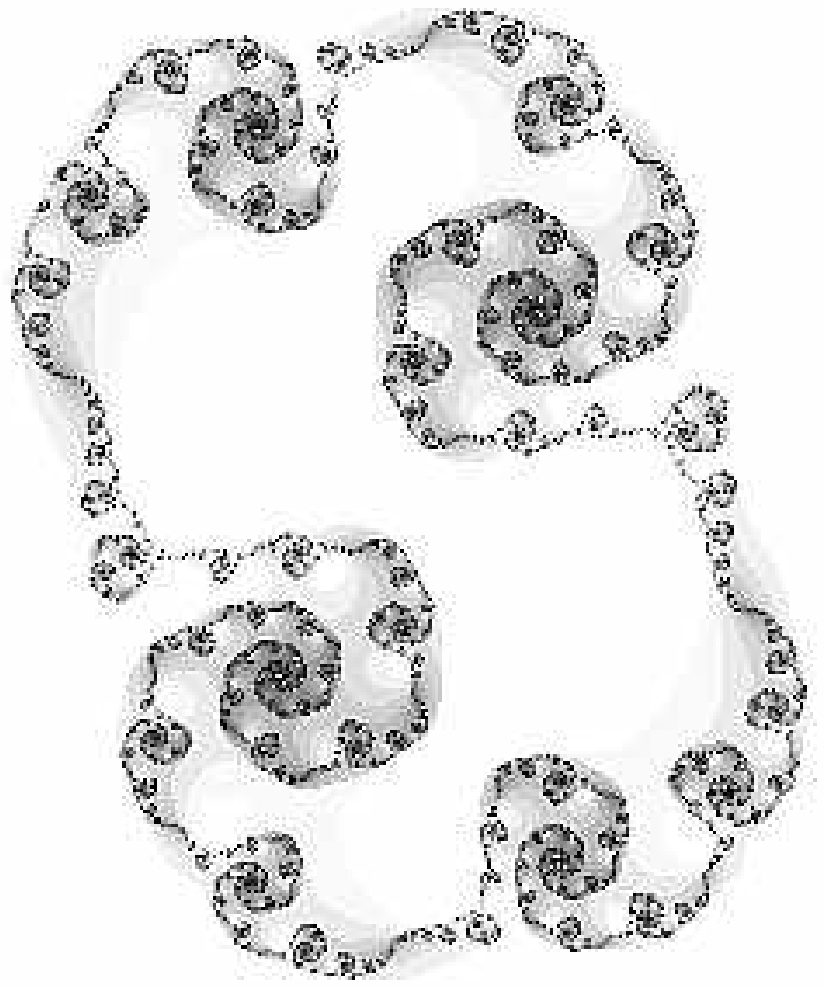}
\caption{The Julia set of $f_a\in V_2$ with $-1\in\d\Omega_0$ and of nearby $f_{a'}\in V_2$
with $-1\in\Omega_0$}
\end{figure}

\subsection{Matings and anti-matings}
Consider two quadratic invariant laminations $L_1$ and $L_2$.
Consider the images $l^{-1}$ of all leaves $l\in L_2$ under the
transformation $z\mapsto 1/z$.
If we straighten all such curves to geodesics in $\{|z|>1\}$,
then we obtain a lamination $L_2^{-1}$ outside the unit disk.
We can form the two-sided lamination $L_1\cup L_2^{-1}$.
The lamination $L_1\cup L_2^{-1}$ is invariant under the map $z\mapsto z^2$
(rather than $z\mapsto 1/z^2$).
This lamination is called the {\em mating} of the laminations $L_1$ and $L_2$.
If the quadratic invariant laminations $L_1$ and $L_2$ correspond
to quadratic polynomials $p_1$ and $p_2$, and if the lamination $L_1\cup L_2^{-1}$
models a rational map $f$, then we say that $f$ is a mating of $p_1$ and $p_2$.
We write $f=p_1\sqcup p_2$ in this case.

This definition of mating is equivalent to the following more standard definition.
Compactify the complex plane by the circle at infinity.
The resulting space is homeomorphic to the closed disk.
Let a polynomial $p_1$ act on one copy of this disk, called $D_1$, and $p_2$
act on another copy, called $D_2$.
Denote by $\gamma_i(t)$ the point on the boundary of $D_i$ of angle $t$.
Identify the boundaries of $D_1$ and $D_2$ by the formula $\gamma_1(t)=\gamma_2(-t)$.
Then the union $D_1\cup D_2$ is homeomorphic to the sphere.
If $p_1$ and $p_2$ have the same degree, then the actions of both polynomials
match on $\d D_1=\d D_2$.
Introduce the minimal equivalence relation $\sim$ on the sphere $D_1\cup D_2$
such that for any point $z\in\d D_1=\d D_2$ that is a common landing point
of two rays, one in $D_1$ and another in $D_2$, the union of these two rays and
the point $z$ belongs to a single equivalence class.
If the quotient $D_1\cup D_2/\sim$ is homeomorphic to the sphere,
and the map $p_1\cup p_2/\sim$ is topologically conjugate to a rational map,
then this rational map is called the {\em mating} of $p_1$ and $p_2$.

Many maps in $V_2$ can be described as matings with the
quadratic polynomial $z\mapsto z^2-1$.
The Julia set of this polynomial is called the {\em basilica}.
The dynamics of $z\mapsto z^2-1$ can be described by a certain
quadratic invariant lamination, which we call the {\em basilica lamination}.
The critical point $0$ of the polynomial $z\mapsto z^2-1$
is periodic of period two: $f(0)=-1$ and $f(-1)=0$.
Thus $z\mapsto z^2-1$ belongs to $V_2$.
Actually, this is the only polynomial of class $V_2$.

\begin{figure}
\centering
\includegraphics[width=6cm]{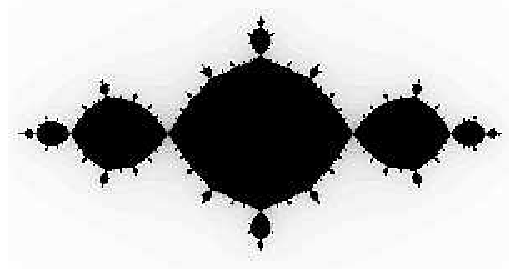}\qquad
\includegraphics[width=6cm]{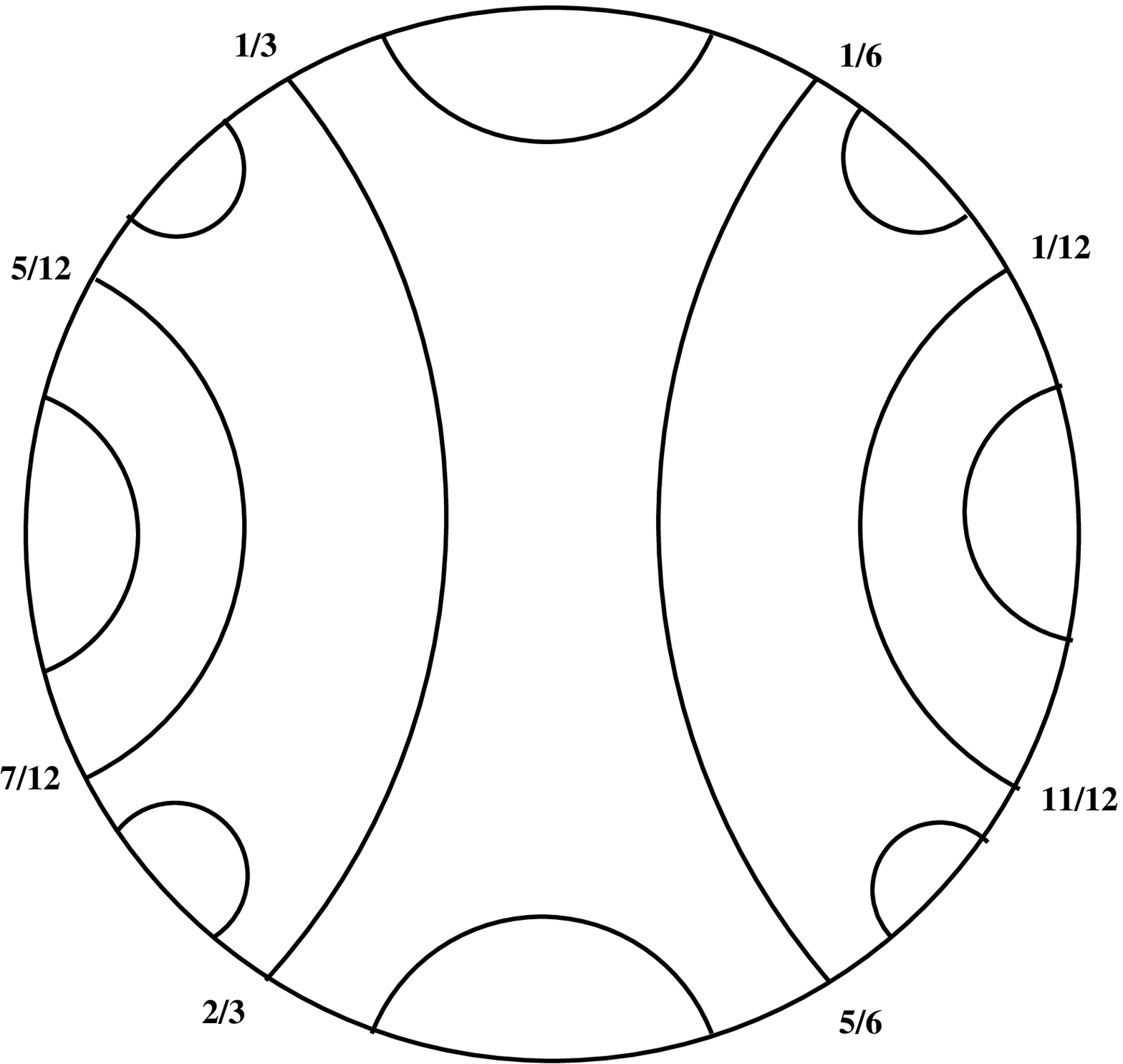}
\caption{The basilica (the Julia set of $z\mapsto z^2-1$) and the basilica lamination}
\end{figure}

\begin{thmB*}
\label{bndry*}
Suppose that the critical point $-1$ of $f_a$
belongs to $\d\Omega_0$ and has angle $\theta_0$.
Let $\theta_0[m]$ denote the $m$-th binary digit of $\theta_0$.
Then, for
$$
y_0=\frac 13\left(1+3\sum_{m=1}^\infty\frac{\theta_0[m]}{4^m}\right),
$$
the mating of the basilica lamination with the lamination $L(-2y_0)$ models
the map $f_a$.
Moreover, the lamination $L(-2y_0)$ itself models a well-defined quadratic
polynomial, so that $f_a$ is a mating of $z\mapsto z^2-1$ with another quadratic polynomial.
\end{thmB*}

The formula for $y_0$ has a simple meaning.
Namely, consider the point in the basilica belonging to the boundary
of the Fatou component of $-1$ and having the (internal) angle $\theta_0$ on this boundary.
Then $y_0$ is the external angle of the same point.
In the terminology we used, internal angles parameterize dynamical rays emanating
from $-1$, whereas external angles parameterize dynamical rays emanating from $\infty$.

Theorem B$^*$ can be deduced from Theorem \ref{bndry}.
Actually, the model with a two-sided lamination invariant under $z\mapsto 1/z^2$
is combinatorially equivalent to the mating model.
However, the model with a two-sided lamination is simpler in some respects.
It can be restated in terms of {\em anti-matings}.
The notion of anti-mating was also introduced by Douady and Hubbard \cite{DH}.
Consider two closed disks $D_1$ and $D_2$ as above,
together with the actions of quadratic polynomials $p_1$ and $p_2$, respectively.
We can glue a topological sphere $D_1\cup D_2$ out of the disks $D_1$ and $D_2$ in the same way as for matings.
However, we define a self-map of this topological sphere differently.
Namely, a point with coordinate $z$ in $D_1$ is mapped to the point with
coordinate $p_1(z)$ in $D_2$.
Similarly, a point with coordinate $w$ in $D_2$ is mapped to the point with
coordinate $p_2(w)$ in $D_1$.
Thus the disks $D_1$ and $D_2$ are interchanged.
Suppose that the quotient of $D_1\cup D_2$ by the equivalence relation $\sim$ is
a topological sphere, and the quotient of the map introduced above is topologically
conjugate to a rational function.
This rational function is then called the {\em anti-mating} of $p_1$ and $p_2$.

Theorem \ref{bndry} can be restated as follows:

\begin{thmB**}
Any map $f_a$ such that $-1\in\d\Omega_0$ is an anti-mating of $z\mapsto z^2$ and
another quadratic polynomial.
\end{thmB**}

From this theorem it follows, in particular, that
the second iteration $f_a^{\circ 2}$ is the mating of two quartic polynomials.
For both of these polynomials, all critical points are either periodic or on the
boundaries of immediate super-attracting basins.

For the case, where the critical point $-1$ is pre-periodic, Theorem \ref{lc} is known,
and the proofs of Theorems \ref{bndry} and B$^*$ are much simpler
(they basically follow from the mating criterion given in \cite{TanLei}).
In this paper, we will concentrate on the case, where $-1$ is not pre-periodic.
As we will see, the angle $\theta_0$ is irrational in this case
(e.g. this follows from Theorem \ref{lc}), however,
we do not assume this a priori.

The results of Theorems \ref{lc}, \ref{bndry} and B$^*$
complement recent results by Aspenberg and Yampolsky \cite{A-Y}.
They prove that any non-renormalizable quadratic polynomial, not in the 1/2-limb and
with all cycles repelling, is mateable with the basilica.
From Theorem B$^*$ it follows, in particular, that any map $f_a$ with $-1\in\d\Omega$
is a mating with a non-renormalizable polynomial, and, therefore, belongs to the class considered in \cite{A-Y}.
The main technical tool of both this paper and \cite{A-Y} are
bubble puzzles suggested by Luo \cite{Luo}.
Luo claimed the main result of \cite{A-Y}, and gave a sketch of a proof, but many
important details were missing.
In other contexts, similar constructions were used in \cite{Yampolsky,Ro}.

The first version of this paper was written before preprint \cite{A-Y} appeared.
It contained a proof of Theorem \ref{bndry} based on a direct construction of the puzzle
specific to our situation.
No analytic continuation was used, but the condition $-1\in\d\Omega$ was essential.
The technique developed in \cite{A-Y} permits to build the puzzle just for
some simple rational maps, and then continue it analytically.
We adopt this approach.

\subsection{The exterior hyperbolic component}
All theorems we stated so far are about maps on the external boundary of $M_2$.
It is natural to attempt studying topology and dynamics of such maps
by approaching them from the {\em exterior component} $\mathcal E$ ---
the only unbounded component of the complement to $M_2$.
There is a simple dynamical description of the set $\mathcal E$:
a map $f_a\in V_2$ belongs to $\mathcal E$ if and only if the free critical point
$-1$ belongs to the immediate basin of the critical cycle $\{0,\infty\}$.
Then we must have $-1\in\Omega_0$, as we will see.

The Julia set of any map $f_a$ in $\mathcal E$ is a quasi-circle, and the restriction
of $f_a$ to the Julia set is conjugate to the map $z\mapsto 1/z^2$.
This follows from a more general theorem of Sullivan \cite{Sullivan1}.
Thus the topology and the dynamics of the Julia set is the simplest possible.
However, a non-trivial combinatorics and a non-trivial dynamics show up
when we consider rays for the second iteration $f_a^{\circ 2}$, and the way
they crash into pre-critical points; more details will come soon.

We give topological models for all maps $f_a$ in $\mathcal E$ in terms of Blaschke products.
The methods used to build these models are not new
(cf. Sullivan and McMullen \cite{McMullen_aut}).
The second iteration $f_a^{\circ 2}$ of the map $f_a$ preserves both components
of the complement to the Julia set.
Pick one particular component.
This is an open topological disk.
Consider a holomorphic uniformization of this topological disk by the round unit disk.
The map corresponding to $f_a^{\circ 2}$ under this uniformization takes the unit disk to itself.
Therefore, it is a quartic Blaschke product.
It is not hard to see that this Blaschke product must actually be the square
of a quadratic Blaschke product
$$
B:z\mapsto z\frac{z+b}{\overline{b}z+1},
$$
where $b$ belongs to the open unit disk.
This gives an idea of how to construct a topological model for $f_a$.

The unit circle divides the Riemann sphere into two disks --- the {\em inside}
and the {\em outside} of the unit circle.
Consider the map $1/B$ that takes the inside to the outside,
and the map $1/z^2$ that takes the outside to the inside.
We would like to glue these maps together but, unfortunately, they do not match on the boundary.
Fortunately, there is a quasi-conformal automorphism $Q$ of the outside of the unit circle
such that the maps $Q\circ 1/B$ and $1/z^2\circ Q^{-1}$ do match on the boundary.
They define a global topological ramified self-covering $g$ of the Riemann sphere of
degree two.
Moreover, there is a natural quasi-conformal structure invariant under $g$.
By the Measurable Riemann Mapping theorem of Ahlfors--Bers \cite{Ahlfors}, the ramified self-covering $g$ is
topologically conjugate to a quadratic rational map.
Clearly, this quadratic rational map must belong to $\mathcal E$.
Conversely, any map in $\mathcal E$ can be obtained by this quasi-conformal surgery.

\subsection{Dynamical rays and external parameter rays}
\label{ss_rays}
Let $f_a$ be a map in $V_2$.
The second iteration $f_a^{\circ 2}$ has two super-attracting fixed points $0$ and $\infty$.
The other four critical points are $-1$, the two preimages of $-1$ under $f_a$,
and $-2$, which is a preimage of $\infty$ under $f_a$.

Consider the Green function $G$ for the map $f_a^{\circ 2}$ that is defined by
the usual formula
$$
G(z)=\lim_{n\to\infty}\frac{\log|f_a^{\circ 2n}(z)|}{2^n}.
$$
This function is negative near 0 and positive near $\infty$.
The gradient of $G$ restricted to the open set $\{G\ne 0\}$ is a smooth vector field that has
singularities at all {\em pre-critical points} (iterated preimages of critical points).
Recall that a {\em ray} is any trajectory of this vector field.

The $\alpha$-limit set of any ray is a single pre-critical point, more precisely,
an iterated preimage of $\infty$ or an iterated preimage of $-1$.
The $\omega$-limit set is either a pre-critical point or a subset of the Julia set
(which is also a single point in a locally connected situation).
If the $\omega$-limit set is a pre-critical point, then this point
is necessarily an iterated preimage of $-1$ (because it can not be an iterated preimage of $\infty$).
Consider any iterated preimage $z$ of $-1$, and assume that $G(z)\ne 0$.
The point $z$ is a saddle point of the Green function.
Thus there are only two rays emanating from $z$ and only two rays crashing into $z$.
The union of the two rays emanating from $z$, together with the point $z$ itself,
is called the {\em ray leaf centered at $z$}.
Thus the ray leaves are in one-to-one correspondence with iterated preimages $z$ of $-1$
such that $G(z)\ne 0$.

Suppose that $a$ belongs to the exterior component $\mathcal E$.
Then the critical point $-1$ of $f_a$ belongs to $\Omega_0$.
Rays emanating from $0$ are parameterized by the {\em angle}.
In a small neighborhood of 0, the map $f_a^{\circ 2}$ is holomorphically conjugate
to the map $z\mapsto z^2$.
Under this local conjugacy, the point 0 is mapped to 0, and
germs of rays are mapped to germs of radial segments.
By definition, the angle of a ray is defined as the angle the corresponding radial
segment makes with the real axis.
We measure angles in radians/$2\pi$.
Thus the measure of the full angle is 1.
Let $R_0(\theta)$ denote the ray of angle $\theta$ emanating from 0.
It is not hard to see that there exists a unique ray $R_0(\theta_0)$ that emanates from $0$ and
crashes into the critical point $-1$.

Fix an angle $\theta_0$.
Consider the set of all parameter values $a$, for which the ray $R_0(\theta_0)$
crashes into the critical point $-1$.
This set is called the {\em external parameter ray of angle $\theta_0$}.
We call an external parameter ray {\em periodic} or {\em non-periodic} according to
whether its angle is periodic or non-periodic under
the doubling map modulo $1$.

M. Rees \cite{Rees} proved that periodic external parameter rays
(except for the zero ray) land at parabolic parameter values.

\begin{bigthm}
\label{param}
All external parameter rays land.
Consider the rational map $f_a\in V_2$
corresponding to the landing point $a$ of a non-periodic external parameter ray of angle $\theta_0$.
For this map, $-1\in\d\Omega_0$.
Moreover, the critical point $-1$ is the point on $\d\Omega_0$ of angle $\theta_0$,
thus the topological dynamics of $f_a$ is described by Theorem \ref{bndry}.
\end{bigthm}

When the first version of this paper was written, I had in mind to deduce this
theorem from Theorem \ref{bndry} by showing that $-1\in\d\Omega$ for all
parameter values on the external boundary, except for countably many parabolic points.
My argument was overly complicated, and I am grateful to M. Lyubich for suggesting
a simpler approach, not using the puzzle.
However, in this paper, Theorem \ref{param} is proved using the parameter
puzzle, a version of that in \cite{A-Y}.
This approach has the advantage that the same combinatorial constructions
are used for both Theorems \ref{lc} and \ref{param}.

\begin{figure}
\centering
\includegraphics[width=5cm]{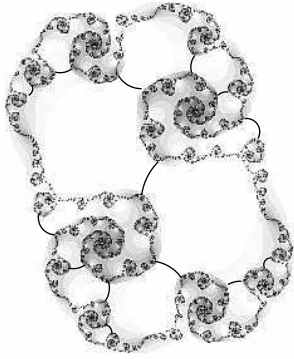}
\caption{Ray leaves for some map in the exterior component of $V_2$}
\end{figure}

\subsection{Ray laminations}
\label{ss_raylam}
Consider a quadratic rational map $f_a$ in the exterior component $\mathcal E$.
Assume that $f_a$ does not lie on a periodic parameter ray
(it can still lie on a strictly pre-periodic parameter ray).
Then each ray leaf of $f_a$ is a curve that is closed in the complement to the Julia set.
The closure of this curve in the Riemann sphere intersects the Julia set in two points ---
the {\em endpoints} of the ray leaf.

Straighten the Julia set to the unit circle, and each ray leaf to a geodesic in the
complement to the unit circle.
Then we obtain a two-sided geodesic lamination.
Since the restriction of the map $f_a$ to the Julia set is conjugate to the map
$z\mapsto 1/z^2$, this two-sided lamination is invariant under $z\mapsto 1/z^2$.
We will call this lamination the {\em ray lamination}.
Ray laminations can be described explicitly.

\begin{bigthm}
\label{ext}
Let $f_a\in V_2$ be a map in the exterior component.
Suppose that $f_a$ lies on a non-periodic external parameter ray of angle $\theta_0$.
Then the ray lamination for $f_a$ coincides with the two-sided lamination $2L(x_0)$,
where
$$
x_0=\sum_{m=1}^\infty\frac{[(2^m-1)\theta_0]+1}{2^{2m+1}}.
$$
\end{bigthm}

We will see that all maps in the same parameter ray give rise to the same ray lamination.
On the other hand, ray laminations corresponding to maps from different
parameter rays can never be the same.

What happens if we approach the external boundary along a non-periodic parameter ray?
The corresponding ray lamination stays the same, but all leaves become shorter and shorter.
In the limit, all leaves of the ray lamination collapse to points.
Thus the same two-sided lamination serves both as a ray lamination
for a map in the exterior component and as a lamination modeling a map on the external boundary.
This picture was the initial motivation for Theorem \ref{bndry} stated above.
However, the formal proof goes differently.
The collapsing of ray leaves can be proved a posteriori, using
Theorems \ref{bndry} and \ref{param}.

\subsection{Hyperbolic components of $V_2$}
From Theorem \ref{param} it follows that the boundary of the exterior component $\Ec$ is a topological circle.
However, it is not a quasi-circle because it has cusps at all parabolic points.
The hyperbolic component $\Ec$ is special because it is the only type II component in $V_2$.
Recall that, according to the terminology of M. Rees \cite{Rees}, a hyperbolic
component in a space of quadratic rational maps is of type II if both critical points
belong to the same cycle.
A hyperbolic component is of type III if one critical point is strictly pre-periodic
and eventually enters the cycle of the other critical point, and of type IV if both
critical points are periodic, with disjoint cycles.
In $V_2$, all type III components are capture components, and all type IV components are
mating components.

Note that the boundaries of type IV components are real analytic curves.
From \cite{A-Y} it follows that the boundaries of type III components are topological circles,
and it is very likely that they are quasi-circles.
Maps on the boundary of a type III component are never critically recurrent.
Thus they exhibit much simpler dynamical behavior, compared with the maps on the external boundary.
On the other hand, maps on the boundary of a type IV component can be much more
complicated, as complicated as quadratic polynomials can be.
In particular, they can have Siegel or Cremer points.

\subsection{A blow-up of $z\mapsto z^2$}
The explicit formula for $x_0$ in terms of $\theta_0$ used in Theorems
\ref{bndry} and \ref{ext} may look mysterious.
We will now explain this formula by describing a simple topological construction it comes from.

Let $z_0$ be any point on the unit circle.
There is a unique probability measure $\mu$ on the unit circle with the
following properties:
\begin{itemize}
\item
The measure $\mu$ is supported on countably many points, namely,
on all iterated preimages of $z_0$ under the map $z\mapsto z^2$
(the point $z_0$ itself is also regarded as an iterated preimage of $z_0$).
\item
For any point $z$ on the unit circle different from $z_0$, we have $\mu\{z^2\}=4\mu\{z\}$.
\end{itemize}
The measure $\mu$ can be given by the following formula
$$
\mu\{z\}=\sum_{m:\ z^{2^m}=z_0}\frac 1{2\cdot 4^m}.
$$
The summation is over all nonnegative integers $m$ such that $z^{2^m}=z_0$.
In particular, if the point $z_0$ is not periodic under the map $z\mapsto z^2$,
then there is at most one summand.
The definition of $\mu$ can be made simple in the non-periodic case:
any preimage of $z_0$ under the map $z\mapsto z^{2^m}$ has measure
$\frac 1{2\cdot 4^m}$.

It is classically known that
there is a unique continuous map $h:S^1\to S^1$ with the following properties:
\begin{itemize}
\item $h(1)=1$, and 1 is in the center of $h^{-1}(1)$.
\item the push-forward of the uniform probability measure under the map $h$ is the measure $\mu$,
\item the map $h$ has topological degree 1.
\end{itemize}
The map $h$ blows up all iterated preimages of the point $z_0$ under $z\mapsto z^2$
in the following sense.
For any point $z$ such that $z^{2^m}=z_0$, the full preimage of $z$ under $h$
is an arc of length $\mu\{z\}$.
In particular, the full preimage $h^{-1}(z_0)$ is a half-circle.
The following proposition is verified by a simple direct computation:

\begin{prop}
\label{blow-up}
If $z_0=e^{2\pi i\theta_0}$ is not periodic under the squaring map $z\mapsto z^2$,
then the half-circle $h^{-1}(z_0)$ is bounded by $e^{2\pi ix_0}$ and $-e^{2\pi ix_0}$, where $x_0$
is expressed in terms of $\theta_0$ by the formula from Theorems \ref{bndry} and \ref{ext}.
\end{prop}

\subsection{Acknowledgements}
I am grateful to M. Lyubich for introducing me to the field of holomorphic dynamics,
for his help and encouragement.
I had stimulating discussions with M. Aspenberg, S. Bonnot, A. Epstein, L. DeMarco,
H. Hakobyan, M. Rees and M. Yampolsky.
M. Aspenberg communicated to me the statements of main results in \cite{A-Y} before they were written.
Finally, I am grateful to J. Milnor and to the anonymous referee for useful remarks and suggestions.

\section{Two-sided laminations $2L(x_0)$}
\label{s_raylam}

{\footnotesize
In this section, we will give details on the explicit
construction of two-sided laminations that appear in Theorems \ref{bndry} and \ref{ext}.
Actually, the construction will be slightly more general, including
the two-sided laminations for parabolic maps, not considered in this paper.}

\subsection{Formulas for $x_0$}
\label{ss_x_0}
Recall that, for a real number $\theta_0$ between 0 and 1 that is not
an odd denominator rational number, we defined the corresponding real number
$x_0$ by the formula
$$
x_0=\sum_{m=1}^\infty\frac{[(2^m-1)\theta_0]+1}{2\cdot 4^m}.
$$
In this subsection, we will find the binary expansion of $x_0$.
Define the functions $\nu_m$ on real numbers between 0 and 1 as follows:
$$
\nu_m(\theta)=\left\{
\begin{array}{cl}
0,& \{2^m\theta\} < \theta\\
1,& \{2^m\theta\} \ge \theta
\end{array}\right.
$$

\begin{prop}
For any real number $\theta$ between 0 and 1, we have
$$
1+[(2^m-1)\theta]=[2^m\theta]+\nu_m(\theta).
$$
\end{prop}

\begin{proof}
There are two cases: $[2^m\theta]=[(2^m-1)\theta]$ and $[2^m\theta]=[(2^m-1)\theta]+1$.
In the first case, subtracting $\theta$ from $2^m\theta$ does not change the
integer part, therefore, $\{2^m\theta\}>\theta$, and $\nu_m(\theta)=1$.
In the second case, subtracting $\theta$ from $2^m\theta$ changes the
integer part, therefore, $\{2^m\theta\}<\theta$, and $\nu_m(\theta)=0$.
\end{proof}

We can now rewrite the formula for $x_0$ as follows:
$$
x_0=\sum_{m=1}^\infty\frac{[2^m\theta_0]}{2^{2m+1}}+
\sum_{m=1}^\infty\frac{\nu_m(\theta_0)}{2^{2m+1}}.
$$
Let us compute the first sum:

\begin{prop}
Let $\theta_0[m]$ denote the $m$-th digit in the binary expansion of $\theta_0$.
Then
$$
\sum_{m=1}^\infty\frac{[2^m\theta_0]}{2^{2m+1}}=\sum_{m=1}^\infty\frac{\theta_0[m]}{2^{2m}}.
$$
\end{prop}

\proof
Denote by $X$ the left hand side of this equality.
Note that the $m$-th binary digit of a real number $\theta$ is equal to $[2^{m}\theta]-2[2^{m-1}\theta]$
for $m\ge 1$.
Therefore, the right hand side is
\[
\sum_{m=1}^\infty\frac{[2^{m}\theta_0]-2[2^{m-1}\theta_0]}{2^{2m}}=2X-X=X.
\eqno{\qed}\]

\medskip

We have proved that
$$
x_0=\sum_{m=1}^\infty\frac{\theta_0[m]}{2^{2m}}+\sum_{m=1}^\infty\frac{\nu_m(\theta_0)}{2^{2m+1}}.
$$
This series represents the binary expansion of $x_0$.
Therefore, we have

\begin{prop}
\label{binary_x_0}
Let $x_0[m]$ denote the $m$-th binary digit of $x_0$.
Then
$$
x_0[2m]=\theta_0[m],\quad x_0[2m+1]=\nu_m(\theta_0)
$$
\end{prop}

\subsection{A forward invariant lamination}
\label{ss_L_0}
Fix a point $z_0=e^{2\pi i\theta_0}$ on the unit circle.
Define a lamination $L_0$ as follows.
We first define a probability measure $\mu$ on the unit circle.
It is given by the following formula:
$$
\mu\{z\}=\sum_{m:\ z^{2^m}=z_0}\frac 1{2\cdot 4^m}.
$$
Next, we consider the map $h$ with the following properties:
\begin{itemize}
\item $h(1)=1$, and 1 is in the center of $h^{-1}(1)$.
\item the push-forward of the uniform probability measure under the map $h$ is the measure $\mu$,
\item the map $h$ has topological degree 1.
\end{itemize}
It blows up all iterated preimages of $z_0$.
We connect two points on the unit circle by a geodesic if these two
points bound the full preimage of a single point under $h$.
The lamination $L_0$ is the set of all such geodesics.
As we will prove shortly, this lamination is {\em forward invariant under $x\mapsto x^4$}: for any
leaf $xy$ of $L_0$, either $x^4=y^4$, or the geodesic $x^4y^4$ is also a leaf of $L_0$.

Note that in the definition of the lamination $L_0$, each leaf $l\in L_0$ comes
together with a specific arc subtended by $l$.
Namely, for a leaf $xy$, the corresponding arc is the full preimage of the point
$h(x)=h(y)$ under the map $h$.
We will call this arc the {\em shadow of the leaf $l$}.
Shadows of different leaves in $L_0$ do not intersect.
Given an arc $\sigma$ on the unit circle, define {\em the bridge over $\sigma$}
as the geodesic connecting the boundary points of the arc $\sigma$.
Thus the bridge over the shadow of a leaf $l\in L_0$ is this leaf $l$ itself.
Denote by $l_0$ the leaf, whose shadow $\sigma_0$ is $h^{-1}(z_0)$.

The lamination $L_0$ has a distinguished gap $G_0$ such that all leaves of
$L_0$ are on the boundary of $G_0$.

\begin{prop}
\label{semi-conj}
The lamination $L_0$ defined above is forward invariant under the map $x\mapsto x^4$.
Moreover, the map $h$ semi-conjugates the endomorphism $x\mapsto x^4$
of the unit circle with the endomorphism $z\mapsto z^2$ everywhere except
on the arc $\sigma_0$.
In other words, $h(x^4)=h(x)^2$ for any point $x$ on the unit circle
such that $h(x)\ne z_0$.
\end{prop}

\begin{proof}
We first define an endomorphism $\phi$ of the unit circle such that $L$ is
forward invariant under $\phi$, and then prove that $\phi$ is the map $x\mapsto x^4$.

Suppose first that a point $x$ on the unit circle does not belong to a
shadow of a leaf of $L_0$.
Then the point $h(x)^2$ has a unique preimage under the map $h$.
Define $\phi(x)$ to be this preimage.
The map $\phi$ thus defined admits a continuous extension that
maps the full $h$-preimage of any point $z$ on the unit circle to the
full $h$-preimage of the point $z^2$, except for $z=z_0$.
To fix one such extension, we require that on each arc
that is the full $h$-preimage of some point, the map $\phi$ act linearly
with respect to the arc-length.
Then $\phi$ is well-defined everywhere except on $\sigma_0$,
and the restriction of $\phi$ to the full $h$-preimage
of any point on the unit circle multiplies all arc lengths by 4.
Indeed, the length of the arc $h^{-1}(z^2)$ is four times bigger that the length
of the arc $h^{-1}(z)$, provided that $z\ne z_0$.
We can also say where $\phi$ should map the arc $\sigma_0$ in order to
be a self-covering of the unit circle.

In the case, where $z_0$ is not periodic under $z\mapsto z^2$,
the arc $\sigma_0$ has length $1/2$.
It should be wrapped twice around the circle under the endomorphism $\phi$.
Both endpoints of $\sigma_0$ should be mapped to the $h$-preimage of $z_0^2$,
which is a single point.
Of course, we require that $\phi$ act linearly on $\sigma_0$.

In the case, where $z_0$ is periodic with the minimal period $p$ under the
map $z\mapsto z^2$, the orbit of the arc $\sigma_0$ under the map $z\mapsto z^4$
consists of $p$ arcs of the following lengths:
$$
\frac4{2(4^p-1)},\ \frac{4^2}{2(4^p-1)},\ \dots,\ \frac{4^p}{2(4^p-1)},
$$
the biggest length being that of $\sigma_0$.
We can arrange that $\sigma_0$ wraps more than twice but less than three times
around the unit disk under the map $\phi$ so that the ends of $\sigma_0$ map to the ends of
the segment of length $4/2(4^p-1)$ (this segment being covered 3 times
by parts of $\sigma_0$ under the map $\phi$).
In all cases, we can arrange that all arc-lengths in $\sigma_0$ get 4 times
bigger modulo $\Z$ under the map $\phi$.

We defined a continuous self-map $\phi$ of the unit circle that is
semi-conjugate to $z\mapsto z^2$ on the complement to the arc $\sigma_0$.
The semi-conjugacy is given by $h$.
It is not hard to see that $\phi$ is a self-covering of the unit circle
and that $\phi(1)=1$.
By definition, the lamination $L_0$ is forward invariant under the map $\phi$.

We will now prove that the map $\phi$ just defined multiplies all arc-lengths
by 4 modulo $\Z$ (in other words, it multiplies all small arc-lengths exactly by 4).
Consider any arc $\sigma$ on the unit circle, whose length is smaller than $1/4$.
We want to show that the length of the arc $\phi(\sigma)$ is 4 times bigger
than the length of the arc $\sigma$.
Since on each arc of the form $h^{-1}(z)$, the map $\phi$ multiplies all arc-lengths by 4,
it suffices to assume that $\sigma$ is the full preimage of the arc $h(\sigma)$ under $h$.
By definition of the measure $\mu$, we have $\mu(h(\sigma)^2)=4\mu(h(\sigma))$.
We also know that $\mu(h(\sigma)^2)$ coincides with the length of the arc $\phi(\sigma)$.
This implies that the length of $\phi(\sigma)$ is $4$ times bigger than the length of $\sigma$.

Since the map $\phi$ multiplies all arc-lengths by 4 and fixes 1, it must have the form
$x\mapsto x^4$.
\end{proof}

\subsection{An invariant lamination}
\label{ss_L}

In this subsection, we extend the lamination $L_0$ to a lamination $L$
invariant under the map $x\mapsto x^4$ in the sense of Thurston.
Recall that a geodesic lamination in the unit disk is said to be {\em invariant}
under the map $x\mapsto x^d$ if
\begin{itemize}
\item it is forward invariant,
\item it is {\em backward invariant}: for any leaf $xy$ of the lamination, there exists
a collection of $d$ disjoint leaves, each connecting a preimage of $x$ with a preimage of $y$
under the map $x\mapsto x^d$.
\item it is {\em gap invariant}:
for any gap $G$, the convex hull $G'$ of the image of $\overline{G}\cap S^1$ is a gap,
or a leaf, or a single point.
\end{itemize}

By a {\em pullback} of a connected set under a continuous map, we mean
a connected component of an iterated preimage of this set.
Recall that the arc $\sigma_0$ was defined as the full preimage of the point
$z_0$ under the map $h$.
The arc $\sigma_0$ is the shadow of some leaf $l_0$.
It is easy to see that the shadow of any other leaf in $L_0$ is a certain
pullback of $\sigma_0$ under the map $x\mapsto x^4$.

\begin{prop}
\label{disjoint_bridges}
Consider the set $A$ of all pullbacks of the arc $\sigma_0$ under the map $x\mapsto x^4$.
The bridges over any two arcs in $A$ are disjoint.
\end{prop}

We need the following lemma:

\begin{lemma}
Consider two different pullbacks $\sigma$ and $\sigma'$ of the arc $\sigma_0$
different from $\sigma_0$.
If the bridges over $\sigma$ and $\sigma'$ intersect, then so do the bridges
over their images under the map $x\mapsto x^4$, unless $\sigma$ or $\sigma'$ coincides
with $\sigma_0$.
\end{lemma}

\begin{proof}
If the bridges over $\sigma$ and $\sigma'$ intersect, then these arcs intersect each other,
but none of them contains the other.
The union $\sigma''$ of the two arcs is also an arc.
If we can show that the length of $\sigma''$ is less than $1/4$, then we
would conclude that the map $z\mapsto z^4$ acts homeomorphically on $\sigma''$,
and hence the images of $\sigma$ and $\sigma'$ have intersecting bridges.

By the {\em depth} of a pullback of $\sigma_0$ we mean the minimal number $n$
such that $\sigma_0$ is the image of the pullback under $x\mapsto x^{4^n}$.
The arcs $\sigma$ and $\sigma'$ cannot be pullbacks of $\sigma_0$ of the same depth,
because different pullbacks of the same depth are disjoint.
By our assumption, neither of the arcs $\sigma$, $\sigma'$ coincides with $\sigma_0$.
Then the length of one arc is at most
$$
\frac 12\left(\frac 14+\frac 1{4^2}+\frac 1{4^3}\dots\right),
$$
while the length of the other arc is at most
$$
\frac 12\left(\frac 1{4^2}+\frac 1{4^3}+\dots\right).
$$
The length of $\sigma''$ is thus at most
$$
\frac 18+\frac 1{4^2}+\frac 1{4^3}+\dots<\frac 14.
$$
This proves the lemma.
\end{proof}

Define the set $A_0$ as the set of all arcs that are shadows of leaves of $L_0$.

\begin{lemma}
\label{A_0_backward_invariant}
The union of the set $A_0$ is backward invariant.
In other words, any pullback of any arc in the set $A_0$ is a subset of some arc in $A_0$.
\end{lemma}

Indeed, this follows from the proof of Proposition \ref{semi-conj}.

\begin{proof}[Proof of Proposition \ref{disjoint_bridges}]
Suppose that there are two arcs from $A$ such that their bridges intersect.
Then, applying to this pair of arcs a suitable iterate of the map $x\mapsto x^4$,
we can make one of the arcs be $\sigma_0$.

Thus we have a pullback $\sigma$ of the arc $\sigma_0$ such that
the bridges over $\sigma_0$ and $\sigma$ intersect.
But this contradicts Lemma \ref{A_0_backward_invariant}.
\end{proof}

We can now define a lamination $L$ as the set of bridges over all pullbacks
of the arc $\sigma_0$.
By Proposition \ref{disjoint_bridges}, the leaves of $L$ are disjoint, so that
$L$ is indeed a lamination.
It is not hard to see that the lamination $L$ does not have any accumulation points
inside the unit disk.

\begin{prop}
The lamination $L$ is invariant under the self-map $x\mapsto x^4$ of the unit circle.
\end{prop}

\begin{proof}
We have already proved the forward and backward invariance.
It remains only to prove the gap invariance.
Define the {\em span $P(l)$ of a leaf $l\in L$} as the open topological disk
bounded by $l$ and the shadow of $l$.
Any gap of $L$ different from $G_0$ can be described as the complement
in a span $P(l)$ to the closures of all spans that lie in $P(l)$.
Denote by $G(l)$ the gap associated with the leaf $l$ in this way.

Suppose that $l$ is a leaf of $L$ different from $l_0$.
Then the image of $l$ under the map $x\mapsto x^4$ is another leaf
$l'$, and the gap $G(l)$ maps to the gap $G(l')$ in the following
sense: the intersection $\overline{G(l)}\cap S^1$ maps to the intersection
$\overline{G(l')}\cap S^1$.
Clearly, the gap $G_0$ maps to itself under the map $x\mapsto x^4$ in this sense.
Moreover, $G_0$ is a {\em critical gap of degree two}: the quotient space $\d G_0/l_0$ maps to
$\d G_0$ as a topological covering of degree two, if we extend the map $x\mapsto x^4$
linearly over leaves.

It remains to consider the gap $G(l_0)$.
This gap is mapped to $G_0$, and this is also a critical gap.
To see that, it is enough to understand what happens with the arc $\sigma_0$,
but this was described in the proof of Proposition \ref{semi-conj}.
\end{proof}

\subsection{A two-sided lamination}
In this subsection, we extend the lamination $L$ to a two-sided lamination $2L$ invariant
under the map $x\mapsto 1/x^2$.
By Proposition \ref{blow-up}, it will be clear that $2L=2L(x_0)$.
In particular, the lamination $2L(x_0)$ exists.

\begin{prop}
\label{antipodal}
The lamination $L$ is invariant under the antipodal map $x\mapsto -x$.
\end{prop}

\begin{proof}
Indeed, if the shadow $\sigma$ of some leaf $l\in L$ is a pullback of
the arc $\sigma_0$ under the map $x\mapsto x^4$,
then $-\sigma$ is also a pullback of $\sigma_0$.
Thus leaves of $L$ map to leaves of $L$ under the map $x\mapsto -x$, and, clearly,
gaps map to gaps.
\end{proof}

Consider the set $L'$ of geodesics outside of the unit circle
connecting pairs of points $1/x^2$ and $1/y^2$, where $x$ and $y$ are endpoints
of a leaf in $L$.

\begin{prop}
The set $L'$ is a geodesic lamination outside of the unit circle.
\end{prop}

Indeed, by Proposition \ref{antipodal}, the images of different leaves in $L$
are either the same or disjoint.

We can now consider the two-sided lamination $2L$ that is the union
of the inside lamination $L$ and the outside lamination $L'$.
By Proposition \ref{blow-up}, we have $2L=2L(x_0)$.

\section{The exterior component}
\label{s_exterior}

{\footnotesize In this section, we describe maps in the exterior component
$\mathcal E$ in terms of a special quasi-conformal surgery performed on Blaschke products.
We also discuss combinatorics of rays.}

\subsection{Anti-matings of Blaschke products}
\label{ss_anti-matings}
Anti-matings of polynomials were considered by Douady and Hubbard \cite{DH}.
In this section, we introduce a similar notion for Blaschke products, together
with an explicit quasi-conformal surgery making these anti-matings into rational functions.
Let $\Delta_0$ denote the inside of the unit circle, and $\Delta_\infty$ the
outside of the unit circle (i.e. the complement to the closed unit disk in the Riemann sphere).
The closures of the open disks $\Delta_0$ and $\Delta_\infty$ are denoted
by $\overline{\Delta}_0$ and $\overline{\Delta}_\infty$, respectively.

A {\em (finite) Blaschke product} is a product of any finite number of
holomorphic automorphisms of the unit disk.
The product here is in the sense of multiplication of complex numbers.
Any holomorphic automorphism of the unit disk extends to a holomorphic
automorphism of the Riemann sphere.
Therefore, Blaschke products are also defined on the whole Riemann sphere.

Consider two Blaschke products $B_0$ and $B_\infty$ of the same degree $d$.
We will make the following assumption on $B_0$ and $B_1$:
{\em the restrictions of these maps to the unit circle are expanding in the
usual metric}.
In particular, this implies that both maps $B_0$ and $B_1$ are hyperbolic.
Let $\alpha_0$ be the restriction of the map $1/B_0$ to the unit circle.
This map takes the unit circle to itself.
Moreover, this is an orientation-reversing self-covering of the unit circle
of degree $-d$ (the negative sign represents the change of orientation).
The restriction $\alpha_\infty$ of the map $1/B_\infty$ to the unit circle
satisfies the same properties.

From a classical theorem of M. Shub \cite{Shub} it follows that
any expanding endomorphism of the unit circle is topologically conjugate to
a map $z\mapsto z^k$; the conjugating homeomorphism is unique (see e.g. \cite{Katok}).
In particular, the maps $\alpha_0$ and $\alpha_\infty$ are topologically
conjugate to the map $z\mapsto z^{-d}$.
Since $\alpha_0$ and $\alpha_\infty$ are $C^\infty$, by \cite{Sullivan},
the conjugating homeomorphism is quasi-symmetric.

The following statement is classical, but we give a proof for completeness:

\begin{lemma}
\label{never_commute}
Consider two endomorphisms of the unit circle, one of which is expanding.
If these two maps have the same topological degree and if they commute, then they coincide.
\end{lemma}

\begin{proof}
The expanding map is conjugate to the map $z\mapsto z^k$ for some $k\ne 0,\pm 1$.
If we lift this map to the universal cover of the unit circle
(i.e. to the real line), then we obtain just the linear map $x\mapsto kx$.
Assume that another map of topological degree $k$ commutes with $z\mapsto z^k$.
The lift of this map to the universal cover has the form $x\mapsto kx+P(x)$,
where $P$ is a periodic function.
Since the two maps commute, we have
$$
(kx)+P(kx)=k(x+P(x)).
$$
Therefore, $kP(x)=P(kx)$, and then $k^nP(x)=P(k^nx)$ for all $n$.
The function $P$ is periodic, hence bounded.
It follows that
$$
P(x)=\lim_{n\to\infty}\frac 1{k^n} P(k^nx)=0
$$
for all $x$.
\end{proof}

Let $\phi$ denote the self-homeomorphism of the unit circle that conjugates
$\alpha_0\circ\alpha_\infty$ with $\alpha_\infty\circ\alpha_0$.
Then we have
$$
\phi\circ\alpha_0\circ\alpha_\infty\circ\phi^{-1}=\alpha_\infty\circ\alpha_0.
$$
From this equation it follows that the maps $\phi\circ\alpha_0$ and
$\alpha_\infty\circ\phi^{-1}$ commute.
By Lemma \ref{never_commute}, this is only possible when
$$
\phi\circ\alpha_0=\alpha_\infty\circ\phi^{-1}.
$$
This is an important functional equation on $\phi$ that we will use.

There is a quasi-conformal self-homeomorphism $Q$ of the disk $\overline{\Delta}_\infty$
that restricts to the map $\phi$ on the unit circle.
This is because $\phi$ is quasi-symmetric: any quasi-symmetric automorphism of
the unit circle extends to a quasi-conformal automorphism of the unit disk, see \cite{Ahlfors}.

Define a self-map $F$ of the unit sphere as follows.
On the disk $\overline{\Delta}_0$, we set $F$ to be $Q\circ (1/B_0)$.
On the disk $\overline{\Delta}_\infty$, we set $F$ to be $(1/B_\infty)\circ Q^{-1}$.
These two maps match on the unit circle by the functional equation on $\phi$.

There is a quasi-conformal structure on the Riemann sphere that is invariant under the map $F$.
Indeed, we can define this structure to be the standard conformal structure
on the unit disk $\Delta_0$, and the push-forward of the standard conformal structure
under $Q$ on the disk $\Delta_\infty$.

By the Measurable Riemann Mapping theorem of Ahlfors and Bers (see \cite{Ahlfors}),
there is a self-homeomorphism of the sphere that takes the quasi-conformal structure we defined
to the standard conformal structure.
Let $f$ be a self-map of the Riemann sphere corresponding to the self-map $F$
under this homeomorphism, and $J$ the image of the unit circle.
The map $f$ is a holomorphic self-map of the Riemann sphere with the Julia set $J$
(which is a quasi-circle).
It has topological degree $d$, hence it is a rational function of degree $d$.

We call the map $f$ the {\em anti-mating} of the Blaschke products $B_0$
and $B_\infty$.

\subsection{The exterior component}
\label{ss_exterior}
In this subsection, we consider one particular example of the general construction
introduced above.
For the map $B_0$, we take a quadratic Blaschke product
$$
B_0(z)=z\frac{z+b}{\overline{b}z+1}
$$
with $|b|<1$.
The origin is a fixed point for this map.
The critical points $c_{1,2}$ of $B_0$ are given by the equation
$
\overline{b}z^2+2z+b=0.
$
Since we have $|c_{1}c_2|=1$, one of the critical points, say $c_1$, satisfies
$|c_1|\le 1$, while for the other critical point $c_2$ we have $|c_2|\ge 1$.
The exact formula for $c_{1,2}$ is
$$
c_{1,2}=\frac{-1\pm\sqrt{1-|b|^2}}{\overline{b}}.
$$
We see that $c_1$ lies in $\Delta_0$, whereas $c_2$ lies in $\Delta_\infty$
(since $|b|<1$, it is clear from this formula that points $c_{1,2}$
cannot both lie on the unit circle).

\begin{prop}
The restriction of $B_0$ to the unit circle is expanding.
\end{prop}

\begin{proof}
By a theorem of Tischler \cite{Tischler}, a Blaschke product $B$ restricts to an expanding
endomorphism of the unit circle if and only if $\lambda B$ has a fixed point
in $\Delta_0$ for all $\lambda$ in the unit circle.
Clearly, the map $B_0$ satisfies this condition.
\end{proof}

For the map $B_\infty$, we just take $z\mapsto z^2$ (the restriction of this map
to the unit circle is obviously expanding).
Let $f=f_{[b]}$ be the anti-mating of the Blaschke products $B_0$ and $B_\infty$.
This is a quadratic rational map.
It depends smoothly (and even real-analytically) on $b$.
However, the dependence is not complex analytic, because the Blaschke product $B_0$
does not depend complex analytically on $b$.

\begin{prop}
The map $f$ has a super-attracting cycle of period two.
\end{prop}

\begin{proof}
Consider the map $F$ from Subsection \ref{ss_anti-matings}.
The image of $0$ under $F$ is $Q(\infty)$, and the image of $Q(\infty)$ is $0$.
Thus $\{0,Q(\infty)\}$ is a periodic cycle of period two for the map $F$.
Moreover, $Q(\infty)$ is a critical point of $F$, hence this cycle is super-attracting.
The map $f$ is quasi-conformally conjugate to $F$.
It follows that $f$ also has a super-attracting cycle of period two.
\end{proof}

This proposition means that $f$ is a map in $V_2$.
In particular, it is holomorphically conjugate to some map of the form
$$
f_a:z\mapsto\frac a{z^2+2z}
$$
or to the map $z\mapsto 1/z^2$.
Thus, for any $b\ne 0$ in the open unit disk, there is a unique complex number
$a$ such that $f_a$ is holomorphically conjugate to $f_{[b]}$.
Recall that $f_{[b]}$ was originally defined only up to a holomorphic conjugacy.
We can fix this degree of freedom by setting $f_{[b]}=f_a$.
For $b=0$, we obtain the map $z\mapsto 1/z^2$.
This defines a map from the unit disk $|b|<1$ to the parameter space $V_2$.
We will call this map the {\em anti-mating parameterization}.
Actually, it is easy to see that each map $f_{[b]}$ belongs to the
exterior component $\mathcal E$ (this is because all critical points of
$f_{[b]}$ are in the immediate basin of attraction of the super-attracting
cycle $\{0,\infty\}$).

\begin{prop}
The anti-mating parameterization is one-to-one:
if maps $f_{[b]}$ and $f_{[b']}$ are holomorphically conjugate, then $b=b'$.
\end{prop}

\begin{proof}
Indeed, if $f_{[b]}$ and $f_{[b']}$ are holomorphically conjugate on the Riemann sphere, then
the squares of the corresponding quadratic Blaschke products
$$
B_0(z)=z\frac{z+b}{\overline{b}z+1},\quad{\rm and}\quad
B'_0(z)=z\frac{z+b'}{\overline{b'}z+1}
$$
are holomorphically conjugate in the unit disk.
Since 0 is the only fixed point for each of the maps $B_0^2$ and ${B'_0}^2$,
a conjugating homeomorphism $\phi$ must fix 0.
Then $\phi$ is just the multiplication by some complex number $\lambda$
such that $|\lambda|=1$.

The point $-b$ is the only preimage of 0 under $B_0$.
Similarly, the point $-b'$ is the only preimage of 0 under $B'_0$.
Therefore, we must have $b'=\lambda b$.
But then the equation $\lambda B^2_0(z)={B'}_0^2(\lambda z)$ yields $\lambda=1$,
after all cancelations.
In particular, $b=b'$.
\end{proof}

We will need the following obvious lemma:

\begin{lemma}
\label{branches}
Let $f$ and $g$ be holomorphic functions defined on some open subsets of $\overline{\C}$.
Suppose that $f$ has no multiple critical points, and fix an open set $U\subseteq\overline{\C}$.
If for every critical value $v$ of $f$, the set $g^{-1}(v)\cap U$ consists of simple critical points of $g$,
then the multi-valued analytic function $f^{-1}\circ g$ has no ramification points in $U$ .
\end{lemma}

In particular, if $U$ is simply connected, $g$ is defined everywhere on $U$,
and $f$ is a ramified covering over $g(U)$, then $f^{-1}\circ g|_U$ splits into single-valued branches.

\begin{prop}
The anti-mating parameterization is onto:
any quadratic rational map of class $\mathcal E$ is holomorphically conjugate
to $f_{[b]}$ for some $b$.
\end{prop}

\begin{proof}
Consider any map $f\in V_2$ in the exterior hyperbolic component $\mathcal E$.
We may assume that $f=f_a$ for some $a$.
Let $\Omega_0$ and $\Omega_\infty$ denote the immediate basins of $0$ and $\infty$,
respectively, for the map $f^{\circ 2}$ (both $0$ and $\infty$ are super-attracting
fixed points for this map).
The proof that $f$ is holomorphically conjugate to (actually, coincides with) some map $f_{[b]}$,
consists of several steps:

{\em Step 1.}
Conjugate $f^{\circ 2}$ by a Riemann map sending $\Omega_0$ to the unit disk
and fixing 0.
The result is a holomorphic self-covering $g$ of the unit disk of degree 4
such that 0 is a fixed critical point and a preimage $-b\ne 0$ of 0 is also
a critical point.
In particular, all preimages of 0 have multiplicity 2, which means by Lemma \ref{branches}
that there is a well-defined holomorphic branch of the function $\sqrt{g}$.
Denote this branch by $B_0$.

{\em Step 2.}
Since $B_0(0)=0$, we conclude that $z\mapsto B_0(z)/z$ is a holomorphic
automorphism of the unit disk that maps $-b$ to 0.
Therefore, it must have the form
$$
\lambda\,\frac{z+b}{\overline{b}z+1},
$$
where $\lambda$ is a complex number such that $|\lambda|=1$.
Conjugating $g$ by a suitable rotation around the origin, we can arrange that
$\lambda=1$ (with a different choice of $b$).

{\em Step 3.}
The map $f^{\circ 2}$ is holomorphically conjugate to $B_0^2$, and hence to $f_{[b]}^{\circ 2}$,
on the set $\Omega_0$.
More precisely, there is a holomorphic embedding $\phi_0:\Omega_0\to\overline{\C}$ such that
$$
\phi_0\circ f^{\circ 2}=f_{[b]}^{\circ 2}\circ\phi_0\quad on\quad \Omega_0.\eqno{(1)}
$$
Moreover, we can assume that $\phi'_0(0)=1$.
In particular, the 0-ray of $f^{\circ 2}$ emanating from $0$ is mapped to the 0-ray
of $f_{[b]}^{\circ 2}$ emanating from 0.
Since the Julia set of $f$ is locally connected, we can extend $\phi_0$ to the closure
of $\Omega_0$.

{\em Step 4.}
All critical values of $f_{[b]}$ are images under $\phi_0$ of the critical values of $f$.
Therefore, by Lemma \ref{branches}, the multi-valued analytic function $f_{[b]}^{-1}\circ\phi_0\circ f$ splits into
two single-valued branches over $\Omega_\infty$.
The $0$-ray for $f$ emanating from $\infty$ gets mapped to the $0$-ray for $f_{[b]}$ emanating from $0$ under $\phi_0\circ f$.
The two preimages of the latter ray under $f_{[b]}$ are the $0$- and $1/2$-rays emanating from $\infty$.
Choose the branch $\phi_{\infty}$ of $f_{[b]}^{-1}\circ\phi_0\circ f$
that takes the $0$-ray emanating from $\infty$ to the $0$-ray for
$f_{[b]}$ emanating from $\infty$.

{\em Step 5.}
The map $\phi_\infty$ is defined on $\Omega_\infty$, and satisfies the following relation:
$$
f_{[b]}\circ\phi_\infty=\phi_0\circ f\quad on\quad\Omega_\infty.\eqno{(2)}
$$
If we substitute this relation into (1), then we obtain the following:
$$
f_{[b]}\circ (\phi_\infty\circ f)=f_{[b]}\circ (f_{[b]}\circ\phi_0)\quad on\quad\Omega_0.
$$
Using the fact that $\phi_\infty$ takes a $0$-ray to a $0$-ray, we conclude that
$$
\phi_\infty\circ f=f_{[b]}\circ\phi_0\quad on\quad\Omega_0.\eqno{(3)}
$$
From formulas (2) and (3) it also follows that
$$
\phi_\infty\circ f^{\circ 2}=f_{[b]}^{\circ 2}\circ\phi_\infty\quad on\quad\Omega_\infty.\eqno{(4)}
$$

{\em Step 6.}
The map $\phi_\infty$ also extends continuously to the Julia set of $f$.
The restrictions of the maps $\phi_0$ and $\phi_\infty$ to the Julia set of $f$
both conjugate the map $f^{\circ 2}$ with $f_{[b]}^{\circ 2}$.
Uniformize both $\Omega_0$ and the basin of $0$ for $f_{[b]}^{\circ 2}$ by the unit disk.
Both maps $f^{\circ 2}$ and $f_{[b]}^{\circ 2}$ correspond to the self-map $z\mapsto z^4$ of the unit circle.
Thus the maps corresponding to $\phi_0$ and $\phi_\infty$ both conjugate $z\mapsto z^4$ with itself.
It follows that these maps differ by a cubic root of unity.
However, both $\phi_0$ and $\phi_\infty$ take the $0$-rays for $f$ emanating from $0$ and $\infty$
to the $0$-rays of $f_{[b]}$ emanating from $0$ and $\infty$, respectively.
Therefore, the restrictions of $\phi_0$ and $\phi_\infty$ to the Julia set of $f$ must coincide.

{\em Step 7.}
We can now define a continuous map
$$
\phi=\left\{\begin{array}{cl}
\phi_0,& on\quad\overline{\Omega}_0\\
\phi_\infty,& on\quad\overline{\Omega}_\infty
\end{array}\right.
$$
By formulas (2) and (3) and their extensions to the Julia set, the map $\phi$ conjugates
$f$ with $f_{[b]}$.
Moreover, $\phi$ is holomorphic on the Fatou set.
It follows that $\phi$ is holomorphic on the Riemann sphere, i.e. $\phi$ is a M\"obius transformation.
Since it fixes $0$, $1$ and $\infty$, the map $\phi$ must be the identity.
We conclude that $f=f_{[b]}$.
\end{proof}

\subsection{Ray dynamics: non-periodic case}
Let $f=f_a$ be a map in the exterior component.
In this subsection, we will study combinatorics of rays for the map $f^{\circ 2}$.

Consider the ray $R_0=R_0(\theta_0)$ in $\Omega_0$ that emanates from $0$ and
crashes into $-1$.
Such ray always exists.
Indeed, there is at least one ray emanating from $0$
that crashes into a pre-critical point (otherwise, the map $f^{\circ 2}$
would be conjugate to the map $z\mapsto z^2$ everywhere on $\Omega_0$).
The pre-critical point this ray crashes into must be an iterated preimage of $-1$.
The image of this ray under the corresponding (necessarily even) iteration of $f$
will be the ray emanating from 0 and crashing into $-1$.

Suppose that the ray $R_0$ is not periodic under the map $f^{\circ 2}$
(i.e. no iterated image of $R_0$ is contained in $R_0$).
This means that the angle $\theta_0$ is not periodic under the doubling.
There are exactly two rays $R_1$ and $R_2$, whose $\alpha$-limit set
is the critical point $-1$.
The images of these rays under the map $f^{\circ 2}$ coincide and
lie on the ray $f^{\circ 2}(R_0)$.

\begin{prop}
\label{R12_land}
The rays $R_1$ and $R_2$ land in the Julia set.
\end{prop}

\begin{proof}
It suffices to prove this for one ray, say, for $R_1$.
First, we need to show that the ray $R_1$ does not crash into pre-critical points.
Assume the contrary: the $\omega$-limit set of $R_1$ is a pre-critical point $x$.
It is an iterated preimage of $-1$, so that we can write $f^{\circ 2n}(x)=-1$
for some positive integer $n$.

The set $f^{\circ 2}(R_1)$ lies on the ray containing $f^{\circ 2}(R_0)$.
Therefore, the set $f^{\circ 2n}(R_1)$ lies on the ray containing $f^{\circ 2n}(R_0)$.
However, the set $f^{\circ 2n}(R_1)$ has the point $-1$ in its closure,
whereas the ray containing $f^{\circ 2n}(R_0)$ does not
(because $R_0$ is not periodic).
A contradiction.

We see that $R_1$ does not crash into pre-critical points.
Therefore, its $\omega$-limit set is a connected subset of the Julia set.
If this subset contains more than one point, then it contains an arc
(i.e. the preimage of an arc under a homeomorphism between the Julia set and the unit circle).
In this case, the $\omega$-limit set of a suitable iterated image of $R_1$
is the whole Julia set.
The iterated images of $R_1$ belong to the rays containing the iterated images of $R_0$.
Thus the $\omega$-limit set of a ray containing a certain iterated image of
$R_0$ is the Julia set.

Consider two strictly pre-periodic rays $R'$ and $R''$ of different
minimal periods emanating from $0$.
If $R_0$ is strictly pre-periodic, we assume additionally that the minimal periods of
$R'$ and $R''$ are different from that of $R_0$.
The rays $R'$ and $R''$ do not crash into pre-critical points, otherwise their suitable iterated
images would belong to the ray $R_0$, which is not pre-periodic or has a different minimal period.
The standard argument of Douady and Hubbard \cite{DH} now applies to show that $R'$ and $R''$ land
in the Julia set (so that their $\omega$-limits are single well-defined points different
from each other).
The closures of the rays $R'$ and $R''$ divide the closed unit disk into two parts,
and the closure of any ray emanating from $0$ can only belong to one part .
This contradicts the statement that the $\omega$-limit set of a certain ray
emanating from $0$ is the whole Julia set.
\end{proof}

\begin{prop}
\label{rays_land}
Any ray for the map $f^{\circ 2}$ either crashes into an iterated
preimage of $-1$ or lands in the Julia set.
\end{prop}

\begin{proof}
Consider any ray $R$.
The $\alpha$-limit set of this ray is an iterated preimage of $0$ or
an iterated preimage of $-1$.
Thus we can map $R$ to a ray emanating from $0$ or from $-1$ by a suitable iteration
of the map $f^{\circ 2}$.
In other terms, we can assume without loss of generality that the ray
$R$ emanates from $0$ or from $-1$.

Consider the first case: $R$ emanates from $0$.
Suppose that $R$ does not crash into a an iterated preimage of $-1$.
Then its $\omega$-limit set is contained in the Julia set.
The rest of the proof goes exactly as in Proposition \ref{R12_land}.
In the second case, the ray $R$ must coincide with $R_1$ or $R_2$.
The result now follows from Proposition \ref{R12_land}.
\end{proof}

Let $\phi$ denote the quasi-symmetric homeomorphism between the unit circle
and the Julia set of $f$ that conjugates the map $x\mapsto 1/x^2$ with the map
$f$:
$$
f(\phi(x))=\phi(1/x^2),\quad x\in S^1
$$
Recall that we defined the two-sided {\em ray lamination} $RL$
associated with $f$ in the following way: $xy\in RL$ if
and only if $\phi(x)$ and $\phi(y)$ are the landing points of
rays emanating from the same iterated $f$-preimage of $-1$.
The geodesic $xy$ is drawn inside or outside of the unit circle
depending on whether this iterated preimage of $-1$ belongs to
$\Omega_0$ or $\Omega_\infty$.

\subsection{Proof of Theorem \ref{ext}}
Consider a parameter value $a$ in the exterior hyperbolic component that does not belong
to a periodic external parameter ray, and the corresponding rational map $f=f_a$.
Let $J$ denote the Julia set of $f$.
We need to prove that the ray lamination $RL$ coincides with some two-sided lamination
$2L(x_0)$ corresponding to a point $z_0=e^{2\pi i\theta_0}$ on the unit circle that is not
periodic under the map $z\mapsto z^2$ (here $x_0$ is expressed through $\theta_0$
as in Theorems \ref{bndry} and \ref{ext}).
To this end, we recover the map $h$ of Subsection \ref{ss_L_0} in terms of $RL$.
We will use the homeomorphism $\phi:S^1\to J$ from the end of the preceding subsection.

For any iterated preimage $z$ of $-1$, we defined the {\em ray leaf}
$Rl(z)$ as the union of $z$ and the two rays emanating from $z$.
Define a continuous map $\tilde h:S^1\to S^1$ as follows:
\begin{itemize}
\item
if $\phi(e^{2\pi i\theta})$ is the landing point of a ray $R_0(\xi)$, then
we set $\tilde h(e^{2\pi i\theta})=e^{2\pi i\xi}$;
\item
otherwise there is a unique ray $R_0(\xi)$ that splits at a precritical point $z$
and such that $Rl(z)\cup J$ separates $0$ from $\phi(e^{2\pi i\theta})$; we set
$\tilde h(e^{2\pi i\theta})=e^{2\pi i\xi}$.
\end{itemize}

\begin{prop}
The map $\tilde h$ coincides with the map $h$ from Subsection \ref{ss_L_0},
with some choice of the point $z_0$.
\end{prop}

\begin{proof}
We will just check that the map $\tilde h$ satisfies all properties of the map $h$.
Since $\phi(1)$ is the landing point of $R_0(0)$, we have $\tilde h(1)=1$.
It is also clear that $\tilde h$ has topological degree 1.
It only remains to verify that the push-forward
of the Lebesgue measure under $\tilde h$ is the measure $\mu$ corresponding to
some point $z_0$ on the unit circle, as it was defined in Subsection \ref{ss_L_0}.
We denote by $\tilde\mu$ the push-forward of the Lebesgue measure under the
map $\tilde h$.

Consider the ray leaf $Rl(-1)=\{-1\}\cup R_1\cup R_2$.
The landing points of rays $R_1$ and $R_2$ divide the Julia set into two arcs.
Choose the arc $\phi(\tilde\sigma_0)$ that is separated from $0$ by $Rl(-1)$.
The arc $\tilde\sigma_0$ of the unit circle has length $1/2$
(because the boundary points of $\phi(\tilde\sigma_0)$ are mapped to the same point under $f$,
and hence the boundary points of $\tilde\sigma_0$ are mapped to the same point under $x\mapsto 1/x^2$).
The image of $\tilde\sigma_0$ under $\tilde h$ is some point $z_0$ on the unit circle
such that $\tilde\mu\{z_0\}=1/2$.
Any ray leaf is an iterated preimage of the leaf $Rl(-1)$.
Therefore, the images under $\tilde h\circ\phi^{-1}$ of all arcs in $J$ subtended by ray leaves are
points on the unit circle that lie in the backward orbit of $z_0$ under the map $z\mapsto z^2$.
Moreover, if $z^{2^m}=z_0$, then we have $\tilde\mu\{z\}=\frac 1{2\cdot 4^m}$.

We see that the measure $\tilde\mu$ coincides with the measure $\mu$ corresponding to the point $z_0$.
Then the map $\tilde h$ is also the same as the map $h$.
\end{proof}

Theorem \ref{ext} follows immediately from this proposition.

\section{Analytic continuation}
\label{s_V2}

{\footnotesize
In this section, we approach the external boundary of $M_2$ from the exterior component.
We will define fixed point portraits for maps on the external boundary using an analytic continuation
argument similar to that in \cite{A-Y}.}

\subsection{The basin of the super-attracting cycle}
Let us first recall the setup.
Our main object is the following family of quadratic rational self-maps of the
Riemann sphere:
$$
f_a(z)=\frac{a}{z^2+2z}.
$$
Infinity is a periodic critical point of period 2 for all maps in this family.
The corresponding orbit is $\{0,\infty\}$.
The other critical point is $-1$.

Denote by $\Omega$ the immediate basin of attraction of the super-attracting cycle $\{0,\infty\}$.
Let $\Omega_0$ and $\Omega_\infty$ be connected components of $\Omega$ containing
0 and $\infty$, respectively.
The restriction of $f_a$ to $\Omega_\infty$ is a 2-fold branched covering of $\Omega_0$.
It follows that $f_a^{-1}(\Omega_0)=\Omega_\infty$.
We will write simply $f$ instead of $f_a$ whenever this notation is unambiguous.
The Julia set of $f$ will be denoted by $J$.

\begin{prop}
\label{-1_not_in_O_infty}
The critical point $-1$ does not belong to the set $\Omega_\infty$.
\end{prop}

\begin{proof}
If $-1\in\Omega_\infty$, then all critical points of $f$ belong to the same Fatou
component.
It is known (see e.g. \cite{Milnor-QuadRat,Rees}) that in this case, the Fatou component containing the
critical points must be invariant, and the Julia set must be totally disconnected.
A contradiction.
\end{proof}

\begin{prop}
\label{top_disks}
Both sets $\Omega_0$ and $\Omega_\infty$ are topological disks.
\end{prop}

\begin{proof}
Consider a small disk $U$ containing the origin.
For any positive integer $n$, define the open set $U_n$ as
the component of $f^{-n}(U)$ containing 0 or infinity depending
on whether $n$ is even or odd.
Since $-1\not\in\Omega_\infty$, each set $U_n$ contains at most one
critical point.
By the Riemann--Hurwitz formula, if $U_n$ is a topological disk,
then $U_{n+1}$ is also a topological disk.
Thus all $U_n$ are simply connected.

The set $\Omega_0$ is the union of $U_n$ for all even $n$.
As the union of a nested sequence of simply connected open sets,
this set is also simply connected.
Similarly, $\Omega_\infty$ is simply connected.
\end{proof}

\subsection{Radial components}
Let $x$ be an iterated preimage of the critical point $\infty$.
It makes sense to talk about {\em rays} emanating from $x$, see Subsection
\ref{ss_rays} for more details.
Every ray hits the Julia set or a pre-critical point (namely, an iterated preimage of
the critical point $-1$).

Define the {\em radial component} of $x$ as the union of $\{x\}$ and all rays emanating from $x$.
We will call the point $x$ the {\em center} of this radial component.
Clearly, every radial component is an open topological disk.
If the critical point $-1$ is not attracted by the cycle $\{0,\infty\}$, then
each radial component is just a Fatou component.
However, the combinatorial structure of radial components is more stable than
that of Fatou components.

Let $A_0$ and $A_\infty$ denote the radial components of $0$ and $\infty$, respectively.
Note that $f(A_\infty)=A_0$, the restriction of $f$ to $A_\infty$ being
a ramified covering of degree 2.
However, in general, the set $f(A_0)$ is strictly contained in $A_\infty$.
The ray of angle $\theta$ emanating from $x$ will be denoted by $R_x(\theta)$.

The following proposition is essentially due to Luo \cite{Luo}:

\begin{prop}
\label{A0_cap_Ainfty}
Suppose that the parameter $a$ is not on the external parameter ray of angle $0$.
Then the intersection of $\overline{A}_0$ and $\overline{A}_\infty$
contains a fixed point $\omega$ of $f$ that is the landing point of
both $R_\infty(0)$ and $R_0(0)$.
\end{prop}

\begin{proof}
First note that if $a$ is not on the external parameter ray of angle $0$,
then the rays $R_\infty(0)$ and $R_0(0)$ both land in the Julia set.
Consider the landing point $\omega$ of the 0-ray in $A_\infty$.
This is a point on the boundary of $A_\infty$ that is
either a fixed point or a point of period 2.
However, the map $f$ has only one orbit of period two, namely, $\{0,\infty\}$.
It follows that $\omega$ is a fixed point.
Since $\omega$ belongs to the boundary of $A_\infty$, it is
also on the boundary of $A_0=f(A_\infty)$.
\end{proof}

Note that the fixed point $\omega$ must be repelling.
Indeed, this fixed point is a univalent function of the parameter defined
on $\C-0$ with the external parameter ray of angle 0 removed.
Since it does not bifurcate over this region, it never becomes parabolic.
Actually, the ramification point for $\omega$ is exactly the puncture $a=0$,
the value of $a$ that does not correspond to any map in $V_2$.

Let $x$ be an iterated preimage of the critical point $\infty$, and $n$
the minimal non-negative integer such that $f^{\circ n}(x)=\infty$.
The number $n$ is called the {\em depth} of $x$ and of the corresponding radial component.
The following statement can be easily deduced from Proposition \ref{A0_cap_Ainfty}
by applying pull-backs under the iterates of $f$:

\begin{prop}
\label{shared_landing}
Suppose that the parameter $a$ is not on the external parameter ray of
a binary rational angle.
Let $A$ be a radial component, and $r\ne 0$ a binary rational angle.
Then the ray of angle $r$ in $A$ lands at a point in the Julia set that is
also the landing point of the $0$-ray in a unique radial component $A'$,
whose depth is bigger than the depth of $A$.
\end{prop}

The uniqueness follows from the following fact:

\begin{prop}
\label{only_one_ray_lands_at_omega}
The ray $R_\infty(0)$ is the only ray in $A_\infty$ landing at $\omega$.
\end{prop}

The proof is similar to that of the following classical statement
about quadratic polynomials: there is only one external ray landing
at the $\beta$ fixed point.

\begin{prop}
If $A$ is a radial component different from $A_\infty$ and $A_0$, then
the fixed point $\omega$ is not in the closure of $A$.
\end{prop}

\begin{proof}
Suppose that $\omega$ is in the closure of $A$.
Then $\omega$ must be the root point of $A$, i.e. the landing point of the zero ray in $A$
(because some ray in $A$ must land at $\omega$, and this can only be the ray of angle zero).
Note that if $A$ has the property $\omega\in\d A$, then $f(A)$ has the same property.
We can now assume that $A$ has the minimal depth among all radial components with this property,
different from $A_\infty$ and $A_0$.
In this case, $A$ must map to $\Omega_\infty$ under the first iteration of $f$,
and the root point of $A$ must coincide with the landing point of $R_\infty(1/2)$.
But this point is different from $\omega$ by Proposition \ref{only_one_ray_lands_at_omega}.
\end{proof}

\begin{cor}
\label{no_share}
Suppose that $-1$ is not an iterated preimage of $\omega$.
Then any iterated preimage of $\omega$ is on the boundary of exactly two radial components.
\end{cor}

This statement can be easily reduced to the preceding proposition by using iterations of $f$.

\subsection{Regulated rays}

Let $r_0,r_1,\dots$ be a finite or infinite sequence of nonzero binary rational angles,
and $x$ an iterated preimage of $\infty$.
Define the set $\Gamma(x,r_0,r_1,\dots)$ as follows.
Let $A_0$ be the radial component centered at $x$.
Start at $x$ and go in $A_0$ along the ray of angle $r_0$ up to the landing point $a_0$.
By Proposition \ref{shared_landing}, the point $a_0$ is the landing point
of the $0$-ray in some radial component $A_1$.
Go along the $0$-ray of $A_1$ to the center of $A_1$.
From the center, go along the ray of angle $r_1$ up to the landing point $a_2$.
Continuing this process (if possible), we obtain a (finite or infinite) sequence of points $a_m$ and radial
components $A_m$ such that $a_m$ is the landing point of the ray
of angle $r_m$ in $A_m$, and, at the same time, the landing point of the $0$-ray in $A_{m+1}$.
We define $\Gamma(x,r_0,r_1,\dots)$ to be the union of the centers of $A_m$,
the rays of angles $r_m$ in $A_m$, the points $a_m$, and the $0$-rays in $A_{m+1}$.
We call $\Gamma(x,r_0,r_1,\dots)$ a {\em regulated ray} starting at $x$.
It is easy to see that there is a continuous embedding $\gamma:[0,\infty)\to\overline{\C}$
such that $\gamma[0,\infty)=\Gamma(x,r_0,r_1,\dots)$ and $\gamma(n+1/2)=a_n$ for all $n=0,1,\dots$.
We say that an infinite regulated ray $\Gamma(x,r_0,r_1,\dots)$ {\em lands} at a point
$z$ if the corresponding path $\gamma(t)$ converges to $z$ as $t\to\infty$.
Note that a regulated ray is well defined unless it crashes into a pre-critical point.
In particular, if the critical point $-1$ is not attracted by the cycle $\{0,\infty\}$,
then all regulated rays are well defined.

\begin{prop}
\label{reg_rays}
Any iterated preimage of $\infty$ can be connected to $0$ or $\infty$ by a
finite regulated ray.
\end{prop}

\begin{proof}
Note that the full preimage of a regulated ray starting at 0 is a pair of regulated rays
starting at $\infty$:
$$
f^{-1}(\Gamma(0,r_1,r_2,\dots))=\Gamma(\infty,r_1/2,r_2,\dots)\cup\Gamma(\infty,(r_1+1)/2,r_2,\dots).
$$
Consider a regulated ray $\Gamma(\infty,r_1,r_2,\dots)$ starting at $\infty$.
The preimage of this ray is the union of $\Gamma(0,r_1,r_2,\dots)$
and a regulated ray starting at $-2$.
But the latter is a part of $\Gamma(\infty,1/2,r_1,r_2,\dots)$.
We see that the preimage of any regulated ray lies in the union of regulated rays.

Using this statement, it is now easy to prove the proposition by induction.
\end{proof}

Note that the intersection of any two regulated rays is an initial segment of both.
The image of a regulated ray starting at $0$ is a regulated ray starting at $\infty$:
$$
f(\Gamma(0,r_1,r_2,\dots))=\Gamma(\infty,r_1,r_2,\dots).
$$
The image of a regulated ray starting at $\infty$ is either a regulated ray
starting at 0 or the union of a regulated ray starting at $\infty$ and the
path between $0$ and $\infty$ along the zero rays of $A_0$ and $A_\infty$.
The latter path will be denoted by $\Gamma[0,\infty]$.
We have
$$
f(\Gamma(\infty,r_1,r_2,\dots))=\left\{
\begin{array}{cl}
\Gamma(0,2r_1,r_2,\dots),& r_1\ne 1/2,\\
\Gamma(\infty,r_2,\dots)\cup\Gamma[0,\infty],& r_1=1/2.
\end{array}
\right.
$$

Let $x$ be the center of some radial component.
The end of a finite regulated ray $\Gamma(x,r_1,\dots,r_n)$
is the center of another radial component, which we will denote by $A(x,r_1,\dots,r_n)$.
By Proposition \ref{reg_rays}, radial components are in one-to-one correspondence with
finite regulated rays starting at $\infty$ or $0$ and such that all angles $r_i$ are nonzero.

\begin{prop}
Let $r_1,r_2,\dots$ be an infinite sequence of binary rational numbers, and suppose that
the parameter $a$ is on the external parameter ray of angle $\theta_0\ne 2^k r_m$.
Then the regulated ray $\Gamma(0,r_1,r_2,\dots)$ is well defined and lands at a point in the Julia set.
\end{prop}

\begin{proof}
The condition $\theta_0\ne 2^kr_m$ guarantees that the regulated ray $\Gamma(0,r_1,r_2,\dots)$
never crashes into a precritical point.
Therefore, it is well defined.
From the hyperbolicity of $f$ it follows that the diameter of $A_m$ decays exponentially,
therefore, the regulated ray lands.
\end{proof}

To emphasize the dependence of a regulated ray on the parameter $a$, we will sometimes
write $\Gamma_a(\infty,r_1,r_2,\dots)$ instead of $\Gamma(\infty,r_1,r_2,\dots)$.
In the sequel, we will need the notion of the {\em angle} of a regulated ray
$\Gamma_a(\infty,r_1,r_2,\dots)$. To define the angle, consider the regulated ray
$\Gamma_1(\infty,r_1,r_2,\dots)$ for the rational map $f_1$, which is M\"obious conjugate to the quadratic
polynomial $p_{-1}:z\mapsto z^2-1$.
The landing point of this ray corresponds to a point in the basilica that
is the landing point of exactly one external ray of angle $\theta$.
We call $\theta$ the angle of $\Gamma(\infty,r_1,r_2,\dots)$.
Clearly, it depends only on the sequence of binary rational numbers $r_1,r_2,\dots$,
not on a specific parameter value $a$.
This definition is parallel to that of \cite{Luo,A-Y}.

\subsection{Fixed point portraits}
For this subsection, the parameter $a$ is in the exterior hyperbolic component, but
not on a rational external parameter ray.

Consider the regulated ray $\Gamma^0=\Gamma(\infty,1/2,1/2,\dots)$.
Note that this regulated ray is contained in its image under $f$.
Therefore, the landing point of it must be a fixed point of $f$.
Denote this point by $\beta$.
For the parameter values under consideration, all periodic points are repelling.
In particular, $\beta$ is a repelling fixed point.

\begin{prop}
\label{beta_ne_omega0}
The fixed point $\beta$ is different from $\omega$.
\end{prop}

\begin{proof}
Suppose that $\beta=\omega$.
Consider a small topological disk $D$ around $\omega$.
We can arrange that the boundary of this disk intersect
each ray $R_\infty(0)$ and $R_0(0)$ at a single point.
Then the union of these rays and $\omega$ divides $D$ into two parts.
The path $\Gamma^0$ lies in one part and is invariant under $f$
(in the sense that $D\cap f(D\cap\Gamma^0)=D\cap\Gamma^0$).
However, the two parts are interchanged under $f$, because the rays
$R_0(0)$ and $R_\infty(0)$ are interchanged.
A contradiction.
\end{proof}

The map $f$ has three fixed points, and we already identified two of them.
Denote the remaining fixed point by $\alpha$ (the notation $\alpha$ and $\beta$ for fixed
points is meant to suggest a similarity with quadratic polynomials).
The $\alpha$-fixed point is the most interesting one.

\begin{prop}
There is a regulated ray landing at the fixed point $\alpha$.
\end{prop}

\begin{proof}
Let $I$ be the closed segment of the ray $R_\infty(\theta_0)$ between the critical value and the
landing point (since $\theta_0$ is irrational, the ray $R_\infty(\theta_0)$ lands in the Julia set).
The map $f^{-1}$ has two well-defined holomorphic branches on the set $\Omega_\infty-I$.

Since $\alpha\ne\omega$, the $\alpha$-fixed point cannot be on the boundary of $A_\infty$.
Consider a ray leaf $Rl$ on the boundary of $A_\infty$ that separates $\infty$
from the fixed point $\alpha$ (this means that any curve in $\Omega_\infty$ connecting
$\infty$ with $\alpha$ must intersect $Rl$).
Let $D$ be the component of the complement to $Rl\cup J$ lying in $\Omega_\infty$
and containing $\alpha$ on its boundary.
There is a holomorphic branch $g$ of $f^{-n}$ mapping $\Omega_\infty-I$ into $D$
and $A_\infty-I$ into a radial component adjacent to $Rl$.
We have $(\Omega_\infty-I)\supset D\supset g(D)\supset g^{\circ 2}(D)\supset\dots$,
and the sets $g^{\circ m}(D)$ converge to $\alpha$ in the Hausdorff metric.

Consider a finite regulated ray $\Gamma(\infty,r_1,r_2,\dots,r_k)$ connecting $\infty$
with the center of the radial component different from $A_\infty$ and
adjacent to the ray leaf $Rl$.
Actually, $k=1$ in our situation (see Proposition \ref{Gamma1}), but this is not important for the time being.
Consider the infinite regulated ray
$$\Gamma_1=\Gamma(\infty,r_1,r_2,\dots,r_k,r_1,r_2,\dots,r_k,\dots),$$
where the sequence of angles is periodic with period $(r_1,r_2,\dots,r_k)$.
Clearly, $g(\Gamma_1)\subset\Gamma_1\cap D$.
Therefore, $\Gamma_1$ lands at the fixed point $\alpha$.
\end{proof}

The map $f^{\circ n}$ takes the path $\Gamma_1$ to itself (modulo the regulated segment $\Gamma[0,\infty]$).
In this sense, $\Gamma_1$ is periodic under $f$.
Denote the minimal period by $q$.
However, $\Gamma_1$ is not fixed, because otherwise $\Gamma_1$
would coincide with the regulated ray $\Gamma(\infty,1/2,1/2,\dots)$ landing at $\beta$.
Consider all images of $\Gamma_1$ under iterations of $f$
(regarded as regulated rays starting at $\infty$ or $0$; the segment $\Gamma[0,\infty]$
appearing in the image should be disregarded), and denote them
by $\Gamma_1$, $\dots$, $\Gamma_q$, where $\Gamma_i=f^{\circ i-1}(\Gamma_1)$.
All regulated rays $\Gamma_i$ land at the fixed point $\alpha$.
The union $\{\alpha\}\cup\Gamma_1\cup\dots\cup\Gamma_q$ is called the
{\em fixed point portrait} for $f$.

With a fixed point portrait consisting of regulated rays
$\Gamma_1,\dots,\Gamma_q$, we associate the set of angles $\{\theta_1,\dots,\theta_q\}$,
where $\theta_i$ is the angle of $\Gamma_i$.

\subsection{Regulated parameter rays}
Let us start with the following landing property:

\begin{prop}
\label{ext_binrat}
Any external parameter ray of a nonzero binary rational angle lands at a parameter $a$,
for which the critical point $-1$ is on the boundary of $\Omega_0$ and
is eventually mapped to the fixed point $\omega$.
\end{prop}

\begin{proof}
Consider an external parameter ray $\Rc$ of a binary rational angle $\theta_0$.
For any parameter $a$ on this ray, the critical value $-a=f_a(-1)$ belongs to the ray $R_\infty(\theta_0)$.
Since $\theta_0$ is strictly pre-periodic under the doubling, the ray $R_\infty(\theta_0)$
lands in the Julia set.
The landing point $z_a$ must be an iterated preimage of the fixed point $\omega$, because there
are no other fixed points on the boundary of $A_0\cup A_\infty$.
The point $z_a$ moves complex analytically (with respect to the parameter $a$)
with finitely many branch points.

Consider any parameter value $a_0$ in the boundary of $\Rc$.
If $a_0$ is not a ramification point for $z_a$, then $z_a$ moves holomorphically
over a neighborhood $O(a_0)$ of $a_0$.
Thus the closure of the ray $R_\infty(\theta_0)$ moves holomorphically
(hence equicontinuously, see \cite{MSS}) over $O(a_0)$.
It follows that, for the parameter value $a_0$, we have $-a_0=z_{a_0}$,
hence it maps eventually to $\omega$.
Clearly, if $a_0$ is a ramification point for $z_a$, then $-1$ is also mapped eventually
to $\omega$.

There are only finitely many parameter values, for which $-1$ is mapped eventually to $\omega$.
It follows that the parameter ray $\Rc$ lands.
For the landing point $a_0$, we must have $-a_0=z_{a_0}$, which can be easily proved
by induction on the exponent of the denominator of $\theta_0$.
The proposition follows.
\end{proof}

Now recall certain facts from \cite{A-Y} that we will use.
We use slightly different language; however, the translation should be straightforward.
The type III hyperbolic components of $V_2$ are in one-to-one correspondence with finite
sequences $(r_1,\dots,r_n)$ of nonzero binary rational numbers.
For each such sequence $(r_1,\dots,r_n)$, the hyperbolic
component $H(r_1,\dots,r_n)$ consists of all parameter values $a$ such that
the critical value $-a$ belongs to the radial component $A(\infty,r_1,\dots,r_n)$.
The dynamical (B\"ottcher) coordinate of $-a$ in $A(\infty,r_1,\dots,r_n)$ defines the
parameter coordinate of $a$ in $H(r_1,\dots,r_n)$.
Thus it makes sense to talk about {\em internal rays} in $H(r_1,\dots,r_n)$:
the internal ray of angle $\theta$ consists of all parameter values
$a\in H(r_1,\dots,r_n)$ such that the critical value $-a$ belongs to the
dynamical ray of angle $\theta$ in $A(\infty,r_1,\dots,r_n)$, or, equivalently,
the ray of angle $\theta$ in a preimage of $A(\infty,r_1,\dots,r_n)$ crashes into
the critical point $-1$.
The parameter value $a$ lies on the boundary of $H(r_1,\dots,r_n)$ if and only if
the corresponding critical value lies on the boundary of $A(\infty,r_1,\dots,r_n)$.
In \cite{A-Y}, this statement is deduced from the $\lambda$-lemma of
Ma\~ne--Sud--Sullivan \cite{MSS}.

Let $\Rc$ be an external parameter ray of a binary rational angle $r$.
Consider the landing point $a$ of $\Rc$.
For the corresponding rational map $f$, the critical point $-1$ lies on the boundary
of $A_0$ and $A_{-2}$, but also on the boundary of $A(\infty,1/2,r)$ and $A(0,r)$.
It follows that $-a$ is on the boundary of $A(\infty,r)$, hence the parameter value $a$
is on the boundary of the type III component $H(r)$.

For a sequence of nonzero binary rational numbers $r_1,r_2,\dots$, define the
{\em regulated parameter ray} $\Delta(\infty,r_1,r_2,\dots)$ as follows.
Start at $\infty$ and go along the external parameter ray of angle $r_1$.
By Proposition \ref{ext_binrat}, this external parameter ray lands at some
point on the external boundary, which is also a boundary point of $H(r_1)$.
Continue along the zero internal ray of $H(r_1)$ up to the center, and then
go along the internal ray of angle $r_2$ up to a boundary point.
It is not hard to see that this boundary point of $H(r_1)$ is also a boundary
point of $H(r_1,r_2)$.
Continue along the zero internal ray in $H(r_1,r_2)$, etc.

The {\em angle} of a regulated parameter ray $\Delta(\infty,r_1,r_2,\dots)$
is defined as the angle of the corresponding regulated dynamical ray
$\Gamma(\infty,r_1,r_2,\dots)$.

\subsection{Analytic continuation of fixed point portraits}
In this subsection, we essentially follow \cite{A-Y}.
Consider a fixed point portrait $\{\alpha\}\cup\Gamma_1\cup\dots\Gamma_q$ with the
set of angles $\{\theta_1,\dots,\theta_q\}$.
The angles $\theta_1,\dots,\theta_q$ divide the unit circle into several arcs.
The shortest complementary arc is called the {\em characteristic arc}.
Suppose that the characteristic arc is bounded by angles $\theta_-$ and $\theta_+$,
taken in the counterclockwise order.
Then it is not hard to see that the critical value $-a$ must lie between
the regulated rays $\Gamma_-$ and $\Gamma_+$ of angles $\theta_-$ and $\theta_+$, respectively.
The following proposition is proved in \cite{A-Y}:

\begin{prop}
\label{param_land}
The regulated parameter rays $\Delta_-$ and $\Delta_+$ of angles $\theta_-$ and
$\theta_+$, respectively, land at a parabolic point not in the closure of
the exterior component.
\end{prop}

The following statement is slightly more general than in \cite{A-Y}, but with similar proof:

\begin{prop}
\label{holom_port}
The fixed point portrait moves holomorphically over the region {\rm (}called a {\em parameter wake)}
bounded by the regulated parameter rays $\Delta_-$ and $\Delta_+$.
\end{prop}

\begin{proof}
For parameter values in the parameter wake, the critical point never enters
a regulated ray of the fixed point portrait.
Therefore, each regulated ray moves holomorphically.
By the $\lambda$-lemma, it follows that the critical portrait also moves holomorphically.
\end{proof}

As a corollary, we have a well-defined fixed point portrait at all points on the external boundary.
Moreover, for any external parameter ray $\Rc$, whose angle is not a binary rational number,
the fixed point portrait moves continuously over $\overline{\Rc}$, and even holomorphically
over some neighborhood of $\overline{\Rc}$.
Note that Proposition \ref{holom_port} fails if we replace regulated rays with
bubble rays (a {\em bubble ray} corresponding to a regulated ray $\Gamma$
is the union of the closures of all radial components intersecting $\Gamma$).
Actually, bubbles (the radial components) do not move continuously on the external boundary.

\subsection{Dynamical and parameter pre-puzzle}
The union of the fixed point portrait and $\Gamma[0,\infty]$ divides the parameter plane
into several pieces, called {\em pre-puzzle pieces of depth 0}.
We use the term pre-puzzle, because we do not employ equipotentials as we should do
to form the actual puzzle pieces.
The point of considering pre-puzzle is that its combinatorics will be stable
along each external parameter ray.
We define {\em pre-puzzle pieces of depth $n$} as $n$-th pull-backs of the pre-puzzle
pieces of depth 0.
By {\em combinatorics} of the pre-puzzle, we mean the information about which rays bound
which pre-puzzle pieces.
The following statement is immediate:

\begin{prop}
The combinatorics of the pre-puzzle stays fixed over each external parameter ray,
whose angle is not a binary rational number.
\end{prop}

Define a parameter pre-puzzle piece of depth $n$ as the locus of parameter values $a$
such that $f_a$ has a given combinatorics of pre-puzzle pieces of depth $\le n$.

\begin{prop}
Every parameter pre-puzzle piece is an open set bounded by
several pairs of regulated parameter rays, each pair having a common landing point.
\end{prop}

\begin{proof}
The proof is straightforward.
Suppose that there is no neighborhood of a parameter value $a_0$, over which
a specified pre-puzzle piece moves holomorphically.
Then, for certain parameter values $a$ in any neighborhood of $a_0$, a certain iterated image
$f_a^{\circ m}(-1)$ of $-1$ enters a regulated ray in the fixed point portrait.
Since a fixed point portrait is invariant, we may assume that $m>0$,
and $f_a^{\circ m-1}(-a)$ lies on some regulated ray in the fixed point portrait.
We conclude that $a_0$ is in the union of closures of finitely many regulated parameter rays.

Thus every parameter pre-puzzle piece is bounded by closures of finitely many
regulated parameter rays.
It is also easy to see that these regulated parameter rays come in pairs, each pair
having a common landing point.
For regulated parameter rays of periodic angles, this follows from Proposition \ref{param_land}.
For regulated parameter rays of strictly pre-periodic angles, this follows from the fact that
the fixed point portrait moves equicontinuously over open neighborhoods of
the closures of such rays.
\end{proof}

As a corollary, we obtain the following

\begin{prop}
\label{same_comb}
Let $\Rc$ be an external parameter ray, whose angle is not binary rational.
For any $n$, there is an open neighborhood $U$ of $\overline{\Rc}$ such that
all pre-puzzle pieces of depth $\le n$ move holomorphically over $U$.
Therefore, for all points in $\overline{\Rc}$, the corresponding rational maps have
the same combinatorics of the pre-puzzle.
\end{prop}

\section{Puzzles, cells and local connectivity}

{\footnotesize
In this section, we deal with maps on the external boundary of $M_2$.
We study two different types of combinatorial partitions for such maps: puzzles and cells.
We need puzzles to prove local connectivity of Julia sets, and cells to establish
topological models.
}

\subsection{Puzzle}
\label{ss_puzzle}
Throughout this section, $f=f_a$ corresponds to a parameter $a$ on the boundary
of some external parameter ray of angle $\theta_0$.
We assume that $\theta_0$ is not binary rational.
Denote by $E_\infty$ some equipotential curve in $A_\infty$ and by
$E_0$ some equipotential curve in $A_0$.
Let $U$ be the component of the complement to $E_\infty\cup E_0$ containing $-1$.
By choosing appropriate equipotentials $E_\infty$ and $E_0$, we
can arrange that $f^{-1}(U)$ be compactly contained in $U$.
Let $\{\alpha\}\cup\Gamma_1\cup\dots\cup\Gamma_q$ be the fixed point portrait for $f$.
{\em Puzzle pieces of depth zero} are defined as connected components
of the complement to the set
$$
G=\Gamma[0,\infty]\cup\bigcup_{i=1}^q \Gamma_i\cup\{\alpha\}\cup E_\infty\cup E_0,
$$
intersecting the Julia set.
A {\em puzzle piece} $P_n$ of any depth $n$ is defined as a connected
component of $f^{-n}(P_0)$, where $P_0$ is a puzzle piece of depth 0.
For any point $z\in J$ not on the boundary of a puzzle piece,
let $P_n(z)$ denote the puzzle piece of depth $n$ containing $z$.
Puzzle pieces $P_n(-1)$ are called {\em critical puzzle pieces}.

A slight variation of this construction leads to the bubble puzzle,
obtained by replacing the regulated rays $\Gamma_i$ with the corresponding bubble rays.
However, we use regulated rays instead of the corresponding bubble rays
because two different bubble rays may touch at iterated preimages of the
critical point $-1$.

\subsection{Rational-like maps}
P. Roesch \cite{Ro_puzzles} generalized the Yoccoz puzzle technique (initially developed for quadratic
polynomials) to a broader class rational maps.
In this subsection, we briefly recall the terminology of \cite{Ro_puzzles}.
Let $U$ and $U'$ be two open sets in $\overline\C$ with smooth boundaries (in particular, both boundaries have finitely
many connected components).
Suppose that $U'$ is compactly contained in $U$.
Consider a proper holomorphic map $f:U'\to U$ with finitely many critical points that
extends to a continuous map from $\overline{U'}$ to $\overline{U}$.
Such a map is called a {\em rational-like map}.
A rational-like map is called {\em simple} if there is exactly one critical point of $f$ in $U'$,
and this critical point is simple.
The filled-in Julia set for $f$ is defined as $\bigcap_{n\ge 0}f^{-n}(U)$.

A finite connected topological graph $G$ is called {\em admissible} for a simple rational-like map $f:U'\to U$
if the following conditions hold:
\begin{itemize}
\item the graph $G$ contains $\d U$ and is contained in $\overline{U}$,
\item the graph $G$ is {\em stable} under $f$, i.e. we have $G\cap U'\subseteq f^{-1}(G)$,
\item the forward orbit of the critical point is disjoint from $G$.
\end{itemize}
A {\em puzzle piece} of depth $n$ (associated with $(G,f,U',U)$) is defined as any
connected component of $f^{-n}(U-G)$.
The collection of all puzzle pieces is called the {\em puzzle}.
This is a generalization of the Yoccoz definition to the case of rational-like maps.
For any point $z$ in $f^{-n}(U-G)$, there is a unique puzzle piece $P_n(z)$ of depth $n$ containing $z$.
If $z$ is a critical point for $f$, then the puzzle pieces $P_n(z)$ are called the
{\em critical puzzle pieces}.

\begin{ex}
Take $f=f_a$, with and $U$ as in the preceding subsection.
Set $U'=f^{-1}(U)$.
Then $f:U'\to U$ is a simple rational-like map.
The topological graph $G$ introduced in the preceding subsection
is admissible for $f$.
\end{ex}

Suppose that the forward orbit of a point $z$ avoids $\Gamma$.
Define the {\em tableau} $T(z)$ of $z$ as the matrix $T(z)_{i,j}=P_i(f^{\circ j}(z))$
of puzzle pieces, where $i$ and $j$ run through all nonnegative integers
(thus the matrix $T(z)$ is infinite down and to the right).
If $z$ is a critical point, then $T(z)$ is called the {\em critical tableau}.
A tableau $T(z)$ is said to be {\em periodic} of period $k$ if
$T_{i,j+k}(z)=T_{i,j}(z)$ for all $i$ and $j$.
The critical tableau $T$ is called {\em recurrent} if the critical point belongs to
$T_{i,j}$ with $j>0$ and $i$ arbitrarily large.

The following theorem is proved in \cite{A-Y}:

\begin{thm}
\label{puzzle_conv}
Let a rational-like map $f:U'\to U$ and a topological graph $G$ be as in Subsection \ref{ss_puzzle}.
If the critical tableau is not periodic, then the critical puzzle pieces $P_n(-1)$ converge to $-1$.
Moreover, for any point $z$, whose forward orbit is disjoint from $G$, the puzzle pieces $P_n(z)$ converge to $z$.
\end{thm}

Below (Subsections \ref{ss_ex} and \ref{ss_annuli}), we sketch a proof of this theorem under the assumption that the parameter
value $a$ belongs to the boundary of an irrational parameter ray (this is what we
actually need for the proof of the main theorems).
For such parameter values, the critical tableau is automatically non-periodic:

\begin{prop}
For the parameter values on the boundary of the external parameter ray of angle $\theta_0$,
the critical tableau is not periodic, provided that $\theta_0$ is not periodic and
not binary rational.
\end{prop}

\begin{proof}
The argument below is similar to that in \cite{Faught}.
Consider a parameter value $a$ on the boundary of the external parameter ray of angle $\theta_0$,
and the corresponding rational map $f=f_a$.
By Proposition \ref{same_comb}, all critical puzzle pieces intersect both $A_0$ and $A_{-2}$.
The intersection of $P_n(-1)$ with $A_0=\Omega_0$ is bounded by two rays in $A_0$ of binary rational angles
$\theta_n^-$ and $\theta_n^+$.
It is easy to see that both $\theta_n^-$ and $\theta_n^+$ converge to $\theta_0$.
Therefore, the intersections of the critical puzzle pieces with $\overline{\Omega}_0$ converge to
the prime end impression of angle $\theta_0$.

From the combinatorics of the puzzle it also follows that the landing points of
binary rational rays in $\Omega_0$ separate the boundary of $\Omega_0$.
In particular, the prime end impressions are disjoint.
If the critical tableau is periodic, then the prime end impression of angle $\theta_0$
for $\Omega_0$ is also periodic.
It follows that $\theta_0$ is periodic, a contradiction.
\end{proof}

An important corollary of this proposition is the following:

\begin{prop}
\label{crit_bndry}
If $\theta_0$ is not periodic, and $a$ is on the boundary of the external parameter
ray of angle $\theta_0$, then the critical point $-1$ lies on the boundary
of the Fatou component $\Omega_0$.
\end{prop}

\subsection{An example}
\label{ss_ex}
Before discussing general combinatorics of puzzles, let us work out one particular example.
We use the same set-up as in Subsection \ref{ss_puzzle}.
Suppose that the regulated rays $\Gamma_i$, $i=1,2,3$, converging to the
fixed point $\alpha$ are
$$
\Gamma_1=\Gamma\left(\infty,\frac 12,\frac 14,\frac 14,\dots\right),\quad
\Gamma_2=\Gamma\left(\infty,\frac 14,\frac 14,\dots\right),\quad
\Gamma_3=\Gamma\left(0,\frac 12,\frac 14,\frac 14,\dots\right).
$$
Consider also preimages of these regulated rays (or, equivalently, regulated rays symmetric to
these regulated rays with respect to $-1$):
$$
\Gamma'_1=\Gamma\left(0,\frac 14,\frac 14,\dots\right),\quad
\Gamma'_2=\Gamma\left(\infty,\frac 34,\frac 14,\frac 14,\dots\right),\quad
\Gamma'_3=\Gamma\left(\infty,\frac 12,\frac 12,\frac 14,\frac 14,\dots\right).
$$

\begin{figure}
\centering
\includegraphics[width=8cm]{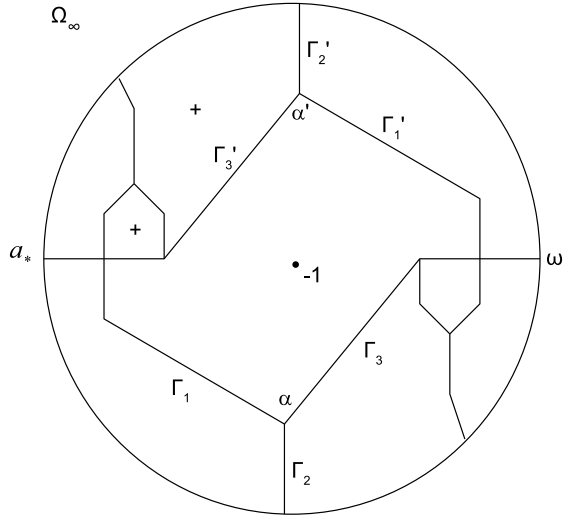}
\caption{An example of the puzzle
(this is a very schematic picture not showing equipotentials and
rays in $\Omega_\infty$)}
\label{pic_puzzle}
\end{figure}

The regulated rays $\Gamma'_1$, $\Gamma'_2$ and $\Gamma'_3$ converge to the point
$\alpha'$ symmetric to $\alpha$ with respect to $-1$, i.e. $\alpha'=-2-\alpha$.
The six paths $\Gamma_i$, $\Gamma'_j$, $i,j=1,2,3$, divide the open
set $U$ into 5 pieces (see Picture \ref{pic_puzzle}).

We see that no puzzle piece of depth 1 is compactly contained in a
puzzle piece of depth 0.
Next, we need to look for puzzle pieces of depth 2 compactly
contained in puzzle pieces of depth 0.
Indeed, there are two puzzle pieces of depth 2 compactly contained in $P^{(0)}(-1)$.
They are marked with sign ``$+$''.

\subsection{Critical annuli}
\label{ss_annuli}
Consider a map $f=f_a$, where the parameter value $a$ is in the closure of an
exterior parameter ray $\Rc$ of an irrational angle $\theta_0$.
Let an open set $U$ be as in Subsection \ref{ss_puzzle}.
In this subsection, we study the rational like map $f:U'\to U$, where
$U'=f^{-1}(U)$, and the puzzle for such map defined in Subsection \ref{ss_puzzle}.

We define the {\em critical annulus of depth $n$} as $R_n(-1)=P_n(-1)-\overline P_{n+1}(-1)$.
If this set is not a topological annulus, we say that the annulus $R_n(-1)$ is {\em degenerate}.

Recall that for quadratic polynomials, the existence of a non-degenerate
critical annulus was settled by the following statement
(see \cite{Milnor,Lyubich}): for a non-renormalizable
quadratic polynomial, the critical orbit enters a non-critical puzzle piece of
depth 1 incident to the point $-\alpha$ (where $\alpha$ is the $\alpha$-fixed point).
There is an analog of this statement for the maps under consideration:

\begin{prop}
\label{depth1}
Let $\alpha'$ be the preimage of $\alpha$ different from $\alpha$, i.e.
$\alpha'=-2-\alpha$.
The critical orbit enters a puzzle piece of depth 1 incident to $\alpha'$
and not containing the critical point $-1$.
\end{prop}

\begin{proof}
Let $\Pi$ be the union of pre-puzzle pieces of depth 1 incident to
$\alpha'$ and not containing the critical point $-1$.
We will write $\Pi_a$ to indicate the dependence of $\Pi$ on the parameter $a$.
We know that the boundary of $\Pi_a$ moves holomorphically with respect to $a$ over
some neighborhood of $\overline\Rc$.

Note that $\Pi_a$ contains either all rays in $A_\infty$
of angles less than $1/2$ or all rays in $A_\infty$
of angles bigger than $1/2$.
Now suppose that $a\in\Rc$.
Then there exists a positive integer $n$ (independent of $a\in\Rc$!) such that
$f_a^{\circ n}(-1)\in\Pi_a$.

Passing to the limit as $a$ approaches the boundary of $\Rc$, we
conclude that $f^{\circ n}(-1)\in\Pi$ for parameter values on the
boundary of $\Rc$.
The proposition follows.
\end{proof}

Unfortunately, unlike the case of quadratic polynomials, not all the puzzle
pieces of depth 1 from Proposition \ref{depth1} are compactly contained in
the critical puzzle piece of depth 0.
Note, however, that the set of angles $2^n\theta_0$ is dense in $\R/\Z$.
In particular, the critical orbit enters all puzzle pieces of depth 1 intersecting $\Omega_\infty$.
Let $\Gamma_1=\Gamma(\infty,r_1,r_2,\dots)$ be a regulated ray landing at $\alpha$.
We can always arrange that $r_1=1/2$ by taking forward images of $\Gamma_1$ under
the iterates of $f$.

For $r_2\ne 1/4, 3/4$, there is a puzzle piece of depth 1 that intersects $\Omega_\infty$
and is compactly contained in the critical puzzle piece of depth $0$.
This is because there are regulated rays in the fixed point portrait of $\alpha$ intersecting
the boundary of $\Omega_\infty$ at points of angles $2^k r_2$, where $k=0,1,2,\dots$.
Since the critical orbit enters this puzzle piece, there is a non-degenerate critical annulus.
We can now use the following theorem \cite{Ro_puzzles} (see also \cite{Milnor-LC} ---
it deals with quadratic polynomials only, but the proof can be taken verbatim in our situation):

\begin{thm}
Suppose that the critical tableau $T$ is recurrent but not periodic.
Also, suppose that there is at least one nondegenerate critical annulus.
Then the critical puzzle pieces converge to the critical point.
\end{thm}

It remains to consider the case, where $r_2$ is $1/4$ or $3/4$ (see Subsection \ref{ss_ex} above).
There are no nondegenerate critical annuli in this case.
Note, however, that the following property holds:

\begin{prop}
The critical point $-1$ can only return to the puzzle piece $P_1(-1)$ under an even iteration of $f$.
\end{prop}

\begin{proof}
It suffices to prove the corresponding statement for the critical pre-puzzle piece of depth 1 and for
parameter values $a$ in the exterior hyperbolic component.
The proposition will follow, if we pass to the limit as $a$ approaches the external boundary of $M_2$.
For the parameter values in the exterior hyperbolic component, the images of $-1$
under all odd iterates of $f$ belong to $A_\infty$, which is disjoint from
the critical pre-puzzle piece of depth 1
(see Picture \ref{pic_puzzle}, in which $\Omega_\infty$ should be replaced with $A_\infty$).
\end{proof}

Therefore, instead of usual critical annuli, we can consider annuli
of the form $P_{n}(-1)-\overline{P_{n+2}(-1)}$, which we call {\em double critical annuli}.
Nondegenerate double critical annuli exist, because there are puzzle pieces of
depth 2 compactly contained in $P_0(-1)$ (see Picture \ref{pic_puzzle}).
We can apply the tableau technique to the double critical annuli.
Namely, the proof of Lemma 1.3 from \cite{Milnor-LC} carries out almost verbatim to double critical annuli.
From this lemma and the Gr\"otzsch inequality, it follows that the critical puzzle pieces
converge to $-1$.

From the convergence of critical puzzle pieces and a simple Koebe distortion argument
it follows that, for any point $z$ in the Julia set of $f$ but not on the boundary
of a puzzle piece, the sequence of puzzle pieces $P_n(z)$ converges to $z$.
The argument goes exactly as for quadratic polynomials.
This concludes the proof of theorem \ref{puzzle_conv}.

\subsection{Cells}
Let $f=f_a$, where $a$ is on the external boundary of $M_2$.
Then, by Proposition \ref{crit_bndry}, the critical point $-1$ is on the boundary of $\Omega$.
In particular, the open set $f^{-1}(\Omega_\infty)$ does not contain critical points.
By the Riemann--Hurwitz theorem, this set consists of two connected components.
One of these components is $\Omega_0$.
The other component contains the point $-2$ (recall that $f(-2)=\infty$).
Denote this component by $\Omega_{-2}$.
Note that in our case, all radial components are Fatou components, e.g.
$A_0=\Omega_0$, $A_\infty=\Omega_\infty$, and $A_{-2}=\Omega_{-2}$.

\begin{prop}
\label{C_*}
The set $\overline{\C}-\overline{\Omega}$ is connected.
\end{prop}

\begin{proof}
Let $C_*$ be the connected component of $\overline{\C}-\overline{\Omega}$ that contains $-2$.
From the existence of the puzzle partition it follows that $C_*$ contains the fixed point $\alpha$.
Indeed, the fixed point portrait contains a regulated ray passing through $-2$ and landing at $\alpha$.
Note also that there is a regulated ray passing through $0$ and landing at $\alpha$.
Therefore, there is a regulated ray passing through $-2$ and landing at $\alpha'=-2-\alpha$
(the point $\alpha'$ is characterized by the properties $f(\alpha')=\alpha$ and $\alpha'\ne\alpha$).
We see that $\alpha'$ also belongs to $C_*$.

The full preimage of $C_*$ under $f$ does not contain critical points.
Therefore, it consists of two connected components.
One of these components contains $\alpha$, and the other component contains $\alpha'$.
It follows that both components are contained in $C_*$, i.e. we have $f^{-1}(C_*)\subset C_*$.

Assume that there is a connected component $V$ of $\overline{\C}-\overline{\Omega}$ different from $C_*$.
The forward orbit of $V$ is disjoint from $C_*$.
Therefore, no iterate of $V$ intersects $G$.
It follows that, for any point $x\in V$ and any depth $n$, we have $V\subset P_n(x)$.
This contradicts the convergence of puzzle pieces, see Theorem \ref{puzzle_conv}.
\end{proof}

The open set $C_*=\overline{\C}-\overline{\Omega}$ is called the {\em main cell}.
Since $-2\in C_*$, we have $\Omega_{-2}\subseteq C_*$.
We define {\em cells of depth $n$} as connected components of $f^{-n}(C_*)$.
Since no cell contains critical points, there are exactly $2^n$ cells of depth $n$.
For any cell $C$ of depth $n$,
there is a unique component of $f^{-n}(\Omega_{-2})$ contained in $C$.
This Fatou component is called the {\em kernel} of the cell.
Note that if a cell has depth $n$, then the depth of its kernel is $n+1$.
Conversely, for each radial component $A$ different from $\Omega_0$ and $\Omega_\infty$,
there is a unique cell containing $A$ as the kernel.
The root point of $A$ (i.e. the landing point of the zero ray in $A$) is also called the {\em root of the cell}.

The methods of Proposition \ref{C_*} also yield the following important result:

\begin{prop}
Any Fatou component of $f$ is eventually mapped to $\Omega_\infty$, i.e. it is a radial component.
\end{prop}

\begin{proof}
Let $V$ be a Fatou component of $f$.
Suppose that the forward orbit of $V$ is disjoint with $\Omega$.
Then, for any point $x\in V$ and any depth $n$, we have $V\subset P_n(x)$.
This contradicts the convergence of puzzle pieces, Theorem \ref{puzzle_conv}.
\end{proof}

We will use cells to encode the dynamics of $f$.
To this end, the following property is crucial:

\begin{thm}
\label{cells}
For any infinite nested sequence of cells
$C^{(1)}\supset C^{(2)}\supset\dots$, the intersection $\bigcap\overline{C^{(n)}}$ consists of a single point.
\end{thm}

We will prove this theorem in Subsection \ref{ss_cells_converge}.
The partition of the Julia set into closures of cells has one major disadvantage:
the critical point $-1$ lies on the boundaries of cells rather than in the interior of a cell.
This is the reason why we also need the puzzle partition.

\subsection{Topology of Fatou components}
In this subsection, we study topology of Fatou components, in particular,
local connectivity and intersection properties of their boundaries.

\begin{prop}
The boundary of $\Omega_\infty$ is locally connected.
\end{prop}

\begin{proof}
We show that all rays in $\Omega_\infty$ land --- the proposition will follow.
It suffices to consider a ray $R$ of irrational angle $\theta$.
Let $I$ be the prime end impression of angle $\theta$ and $x$ any point of $I$.
Clearly, the forward orbit of $I$ is disjoint from $G$.
It follows that, for any depth $n$, we have $I\subset P_n(x)$.
By Theorem \ref{puzzle_conv}, the puzzle pieces $P_n(x)$ converge to $x$.
It follows that $I=\{x\}$, and that $R$ lands at $x$.
\end{proof}

\begin{prop}
\label{no_common_landing}
Two different rays in $\Omega_\infty$ cannot land at the same point.
\end{prop}

\begin{proof}
Assume the contrary: there are two rays in $\Omega_\infty$ that land at the same point.
The union of these rays, the common landing point and $\infty$ divides the Riemann sphere
into two parts.
Each part must contain points of the complement to $\overline{\Omega}_\infty$.
This contradicts Proposition \ref{C_*}.
\end{proof}

\begin{prop}
\label{only_-1}
The critical point $-1$ is the only intersection point of $\overline{\Omega}_0$
and $\overline{\Omega}_{-2}$.
\end{prop}

\begin{proof}
By our assumption, the critical point $-1$ belongs to the boundary of $\Omega_0$.
Note that the map $z\mapsto -2-z$ takes $\Omega_0$ to $\Omega_{-2}$.
It follows that $-1$ is on the boundary of $\Omega_{-2}$, therefore,
$-1\in\overline{\Omega}_0\cap\overline{\Omega}_{-2}$.

Suppose that $x_0\ne -1$ is another point in $\overline{\Omega}_0\cap\overline{\Omega}_{-2}$.
Let $R_0$ and $R_{-2}$ be rays in $\Omega_0$ and $\Omega_{-2}$, respectively, that land at $x_0$.
If $R_0$ and $R_{-2}$ map to the same ray under $f$, then $x_0$ must be a critical point
(indeed, $f$ is not injective in any neighborhood of $x_0$).
Suppose that $f(R_0)$ and $f(R_{-2})$ are different rays in $\Omega_\infty$.
However, they land at the same point $f(x_0)$.
Contradiction with Proposition \ref{no_common_landing}.
\end{proof}

\begin{prop}
\label{only_omega}
The fixed point $\omega$ is the only intersection point of $\overline\Omega_0$ and
$\overline\Omega_\infty$.
\end{prop}

\begin{proof}
Assume the contrary: $x$ is another point in $\overline\Omega_0\cap\overline\Omega_\infty$.
The union of $\{0,x,\omega,\infty\}$ and the rays in $\Omega_0$ and $\Omega_\infty$
landing at $\omega$ and $x$ is a simple closed curve.
This curve divides the Riemann sphere into two parts.
By Proposition \ref{C_*}, only one part can contain points of $C_*$.
Then, in the other part, the boundaries of $\Omega_0$ and $\Omega_\infty$ coincide.
It is easy to see that, in this case, $\d\Omega_0=\d\Omega_\infty$, a contradiction.
\end{proof}

\subsection{Topology of cells}
There are two cells of depth 1.
Denote them by $C_0$ and $C_1$.
Let $a_*$ be the landing point of the ray $R_\infty(1/2)$.
This point belongs to the boundary of both $\Omega_\infty$ and $\Omega_{-2}$.
The following is a consequence of Propositions \ref{only_-1} and \ref{only_omega}.

\begin{prop}
The intersection of $\overline C_0$ and $\overline C_1$ is $\{\omega,a_*,-1\}$.
\end{prop}

For the following, we need two simple lemmas.

\begin{lemma}
\label{kernel}
The kernel of any cell $C$ is the Fatou component of the biggest depth contained in $C$.
\end{lemma}

\begin{proof}
Any cell gets eventually mapped to the main cell $C_*$.
For the main cell, the statement is obvious.
\end{proof}

\begin{lemma}
\label{depth}
Let $C$ be a cell.
The depth of any Fatou component lying in $C$ is bigger than that of $C$.
\end{lemma}

\begin{proof}
Suppose that $C$ has depth $n$.
Then the kernel of $C$ has depth $n+1$.
The statement now follows from Lemma \ref{kernel}.
\end{proof}

We need to establish convergence of certain nested sequences of puzzle pieces.
For any positive integer $n$, there are two puzzle pieces of depth $n$ containing
the fixed point $\omega$ on their boundary.
One of these puzzle pieces, say, $P_{n,0}(\omega)$, intersects $C_0$, and
the other puzzle piece, $P_{n,1}(\omega)$, intersects $C_1$.
We have
$$
P_{n+1,0}(\omega)\subset P_{n,0}(\omega),\quad P_{n+1,1}(\omega)\subset P_{n,1}(\omega).
$$

\begin{prop}
\label{puzzle_omega}
The nested sequence of closed sets $\overline P_{n,0}(\omega)$ converges to $\omega$.
Similarly, the nested sequence $\overline P_{n,1}(\omega)$ converges to $\omega$.
\end{prop}

\begin{proof}
Let $x$ be any point in $\bigcap_{n\ge 1}\overline P_{n,0}(\omega)$.
If $x$ is different from $\omega$, then it is easy to see that the forward
orbit of $x$ is disjoint from $G$.
By Theorem \ref{puzzle_conv}, it follows that the puzzle pieces $P_n(x)$
converge to $x$.
However, we must have $P_n(x)=P_{n,0}(\omega)$.
This is a contradiction, which shows that the sequence $\overline P_{n,0}(\omega)$ converges to $\omega$.
A proof that $\overline P_{n,1}(\omega)$ converges to $\omega$ is similar.
\end{proof}

For any positive integer $n$, there are exactly 2 cells of depth $n$
that contain the critical point $\omega$ on the boundary.
One of these cells, say, $C^{(n)}_0(\omega)$, is contained in $C_0$,
and the other cell, $C^{(n)}_1(\omega)$ is contained in $C_1$.
It is easy to see that, for $n\ge 1$, the kernel of every cell $C^{(n)}_0(\omega)$
and $C^{(n)}_1(\omega)$ touches both $\Omega_0$ and $\Omega_\infty$
(we say that two open sets {\em touch} if their closures intersect).
Moreover, either the point where the kernel of $C^{(n)}_0(\omega)$ touches $\Omega_0$
or the point where the kernel of $C^{(n)}_0(\omega)$ touches $\Omega_\infty$ is eventually mapped to $\omega$
(the same holds for $C^{(n)}_1(\omega)$).
In particular, there is a finite regulated ray $\Gamma(0,r_1)$ or
$\Gamma(\infty,r_1)$ passing through $\Omega_0$ or $\Omega_\infty$ and the
kernel of the cell $C^{(n)}_0(\omega)$.

\begin{prop}
\label{omega}
The nested sequence of closed sets $\overline C^{(n)}_0(\omega)$ converges to $\omega$.
Similarly, the nested sequence $\overline C^{(n)}_1(\omega)$ converges to $\omega$.
\end{prop}

\begin{proof}
We show that every puzzle piece $P_{n,0}(\omega)$ contains some cell $C^{(m)}_0(\omega)$.
The proposition will follow then from Proposition \ref{puzzle_omega}.
The puzzle piece $P_{n,0}(\omega)$ is bounded by a finite number of
regulated rays and equipotentials.
We can choose $m$ such that the kernel $A$ of the cell $C^{(m-1)}_0(\omega)$
touches both $\Omega_0$ and $\Omega_\infty$ at interior points of $P_{n,0}(\omega)$.
Suppose that, say, $\Gamma(0,r_1)$, is a finite regulated ray passing though
$\Omega_0$ and the kernel of the cell $C^{(m-1)}_0(\omega)$
(if it does not exist, then there is a regulated ray $\Gamma(\infty,r_1)$ passing
though $\Omega_\infty$ and the kernel of the cell $C^{(m-1)}_0(\omega)$).
Since the intersection of any pair of regulated rays is an initial segment of both,
the regulated ray $\Gamma(0,r_1)$ is disjoint from the boundary of the puzzle piece $P_{n,0}(\omega)$.
Therefore, the kernel of the cell $C^{(m-1)}_0(\omega)$ is disjoint from the boundary of
the puzzle piece $P_{n,0}(\omega)$, and the cell $C^{(m)}(\omega)$ is contained in $P_{n,0}(\omega)$.
\end{proof}

\subsection{Cells converging to $\alpha$}
For any point $x$ in the Julia set of $f$ but not on the boundary of a cell, there
is a unique cell $C^{(n)}(x)$ of depth $n$ containing $x$.

\begin{prop}
\label{alpha}
The nested sequence of cells $C^{(n)}(\alpha)$ containing $\alpha$ converges to $\alpha$, i.e.
$$
\bigcap_{n=1}^\infty\overline{C^{(n)}(\alpha)}=\{\alpha\}.
$$
\end{prop}

\begin{proof}
The proof consists of several steps.

{\em Step 1.} Let $A_n$ denote the kernel of $C^{(n)}(\alpha)$.
If all $A_n$ touch $\Omega_0$ and $\Omega_\infty$, then $C^{(n)}(\alpha)$
coincide with $C_0^{(n)}(\omega)$ or with $C_1^{(n)}(\omega)$.
However, this contradicts Proposition \ref{omega}.

{\em Step 2.} It follows that some $A_n$ does not touch $\Omega_0$ or
does not touch $\Omega_\infty$.
Suppose that $n$ is the minimal index with this property.
Then it is easy to see that $A_n$ touches $\Omega_\infty$ and $\Omega_{-2}$.
It follows that $C^{(n)}(\alpha)$ does not touch $\Omega_0$.

{\em Step 3.}
Consider the intersection $I$ of $\overline C^{(n)}(\alpha)$.
This is a compact connected subset of the Julia set for $f$.
By step 2, the set $I$ is disjoint from $\overline\Omega_0$.
Since $I$ is forward invariant under $f$, it is also disjoint
from $\overline\Omega_\infty$.
By the same reason, $I$ is disjoint from $\overline\Omega_{-2}$.

{\em Step 4.}
It follows that there is a cell $C^{(n)}(\alpha)$ that does not
touch $\Omega_\infty\cup\Omega_0\cup\Omega_{-2}$.
There is a single valued branch of $f^{-n+1}$ that takes
the cell $C^{(1)}(\alpha)$ (which is $C_0$ or $C_1$) to $C^{(n)}(\alpha)$.
Note that $C^{(n)}(\alpha)$ is compactly contained in $C^{(1)}(\alpha)$.
The proposition now follows from the Poincar\'e distance argument.
\end{proof}

\begin{prop}
\label{Gamma1}
There is a regulated ray of the form
$$
\Gamma(\infty,1/2,r_2,r_2,\dots,r_2,\dots)
$$
converging to $\alpha$.
\end{prop}

\begin{proof}
Let $a_n$ be the root point of the cell $C^{(n+1)}(\alpha)$.
By Proposition \ref{alpha}, not all points $a_n$ are on the boundary of $\Omega$,
whereas $a_{-1}=a_*$ is in the boundary of $\Omega_\infty$.
It follows that there is a nonnegative integer $n$ such that $a_n\in\d\Omega_{-2}$.
Let $r_2$ be the angle of $a_n$ with respect to $\Omega_{-2}$.
Consider the regulated ray $\Gamma_1=\Gamma(\infty,1/2,r_2,r_2,\dots)$.
This ray is periodic; let $q$ be the minimal period.
There is a branch $g$ of $f^{-q}$ that takes $C^{(n+1)}(\alpha)$ to
$C^{(n+q+1)}(\alpha)\subset C^{(n+1)}(\alpha)$.
Clearly, we have $g(\Gamma_1)\subset\Gamma_1$.
The proposition now follows.
\end{proof}

\section{Topological models}
\label{s_topmod}

{\footnotesize In this section, we give topological models for rational maps $f=f_a$
satisfying the condition $-1\in\d\Omega$.  We use the partition
of the Julia set into cells to encode the topological dynamics of $f$.}

\subsection{Convergence of cells and a proof of Theorem \ref{lc}}
\label{ss_cells_converge}

In this section, we prove Theorem \ref{cells}: all nested sequences of cells converge to singletons.

\begin{prop}
Consider any point $z$ in the Julia set of $f$ different from $\alpha$
and such that the forward orbit of $z$ is disjoint from $\{-1,\omega\}$.
Then there is a cell $C(z)$ that contains $z$ in its closure and lies in a puzzle piece of depth $0$.
\end{prop}

\begin{proof}
Since $z$ does not coincide with $\alpha$, it avoids the closure of a cell
$C^{(n)}(\alpha)$ containing $\alpha$ (this follows from Lemma \ref{alpha}).
Let $N$ denote the maximal depth of a Fatou component intersecting
some regulated ray $\Gamma_i$ but not lying in the cell $C^{(n)}(\alpha)$.
It is not hard to see that the cell $C(z)=C^{(N)}(z)$ of depth $N$ lies
in some puzzle piece of depth 0, see Lemma \ref{depth}.
By definition, $z$ belongs to $\overline{C^{(N)}(z)}$.
\end{proof}

The following statement now follows from the convergence of puzzle pieces.

\begin{prop}
\label{conv}
Let $z$ be any point in the Julia set of $f$, whose forward orbit is disjoint from
$\{-1,\omega,\alpha\}$.
We have
$$
\bigcap_{n=1}^\infty\overline{C^{(n)}(z)}=\{z\}.
$$
\end{prop}

Note that iterated preimages of $\omega$ are the only points in the Julia set
that lie on the boundaries of puzzle pieces.

Let $z$ be an iterated preimage of $-1$.
Then, for each depth $n$, there are two cells $C^{(n)}_0(z)$ and $C^{(n)}_1(z)$
having $z$ on the boundary.
We can arrange the indexing so that to have
$$
C^{(n+1)}_0(z)\subset C^{(n)}_0(z),\quad C^{(n+1)}_1(z)\subset C^{(n)}_1(z).
$$
We will also assume that
$$
C_0^{(n)}(-1)\subseteq C_0,\quad C_1^{(n)}(-1)\subseteq C_1.
$$

\begin{prop}
\label{precrit}
For any iterated preimage $z$ of the critical point $-1$, we have
$$
\bigcap_{n=1}^\infty \overline{C^{(n)}_0(z)}=\bigcap_{n=1}^\infty \overline{C^{(n)}_1(z)}=\{z\}.
$$
\end{prop}

\begin{proof}
It suffices to prove this for $z=-1$.
Note that $C_0^{(n)}(-1)$ and $C^{(n)}_1(-1)$ are centrally symmetric with respect to $-1$.
If, say, $\alpha\in C_0$, then $C_1$ is contained in a single puzzle piece of depth $0$,
namely, in the critical puzzle piece $P_0(-1)$.
The critical orbit returns to $\overline{C}_1$, and hence to $P_0(-1)$, infinitely many times.
Suppose that $f^{\circ m}(-1)\in \overline{C}_1$.
Then, by the pullback argument, $C^{(m)}_0(-1)$ or $C^{(m)}_1(-1)$
is contained in $P_{m-1}(-1)$,
which is the pullback of $P_0(-1)$ along the critical orbit.
Since $m$ can be made arbitrarily large, the diameters of $C^{(n)}_0(-1)$ and
$C^{(n)}_1(-1)$ tend to $0$ as $n\to\infty$.
\end{proof}

\begin{proof}[Proof of Theorem \ref{cells}]
Consider a nested sequence of cells $C^{(n)}$.
The intersection of all $\overline{C^{(n)}}$ is non-empty.
Let $z$ be any point in this intersection.
If $z$ is not in the backward orbit of $\{-1,\omega,\alpha\}$, then the convergence follows
from Proposition \ref{conv}.
If $z$ is an iterated preimage of $\omega$, then the convergence follows
from Proposition \ref{omega}.
If $z$ is an iterated preimage of $-1$, then the convergence follows
from Proposition \ref{precrit}.
Finally, if $z$ is an iterated preimage of $\alpha$, then the convergence follows
from Proposition \ref{alpha}.
\end{proof}

Note that Theorem \ref{lc} follows from Theorem \ref{cells}, because
the intersection of the Julia set with each cell is connected.
Indeed, any component of the complement to this intersection is a simply connected Fatou component.

\subsection{Encoding of the Julia set}
\label{ss_encoding}
In this subsection, we encode all points of the Julia set by binary sequences.
Our main tool is Theorem \ref{cells}.
Consider a cell $C$ of depth $n$.
The {\em address of $C$} is a finite binary sequence $\eps_1\dots\eps_n$ defined as follows.
We set $\eps_k=0$ or 1 depending on whether $f^{\circ k-1}(C)$ is contained in $C_0$ or in $C_1$.
We will think of the main cell as having the empty address.
For any finite binary sequence $\eps_1\dots\eps_n$, there is a unique cell
$C_{\eps_1\dots\eps_n}$ with address $\eps_1\dots\eps_n$.
We have $f(C_{\eps_1\eps_2\dots\eps_n})=C_{\eps_2\dots\eps_n}$.

We can now define a continuous map from all infinite binary sequences
to the Julia set of $f$ (the set of infinite binary sequences is considered
as a topological space with respect to the direct product topology).
Given an infinite binary sequence $\eps_1\dots\eps_n\dots$, define
the point $z_{\eps_1\dots\eps_n\dots}$ to be the only point in
$\bigcap_{n=1}^\infty \overline{C}_{\eps_1\dots\eps_n}$.
We have
$$
f(z_{\eps_1\eps_2\dots\eps_n\dots})=z_{\eps_2\dots\eps_n\dots}.
$$
The sequence $\eps_1\dots\eps_n\dots$ is called an {\em address}
of the point $z_{\eps_1\dots\eps_n\dots}$.
Note that the same point can have different addresses.

From now on, we will assume that the cells $C_0$ and $C_1$ of depth 1
are indexed so that the landing points of all rays $R_\infty(\theta)$
with $\theta<1/2$ belong to the closure of $C_0$.
Then the landing points of all rays $R_\infty(\theta)$ with $\theta>1/2$
belong to the closure of $C_1$.
Clearly, this can be arranged.

\begin{prop}
\label{address_of_-1}
The critical point $-1$ is encoded by exactly two binary sequences, namely,
$$
-1=z_{0\eps^*_1\dots\eps^*_n\dots}=z_{1\eps^*_1\dots\eps^*_n\dots},\quad
\eps^*_{2m}=\theta_0[m],\quad \eps^*_{2m+1}=1-\nu_m(\theta_0),
$$
where $\theta_0[m]$ denotes the $m$-th digit in the binary expansion of $\theta_0$,
and the function $\nu_m$ is that introduced in Subsection \ref{ss_x_0}.
\end{prop}

\begin{proof}
The point $-1$ belongs to the closures of both $C_0$ and $C_1$.
However, the remaining address of $-1$ is well-defined:
the $m$-th digit is $0$ if $f^{\circ m-1}(-1)$ belongs to $\overline{C}_0$ and
$1$ if $f^{\circ m-1}(-1)$ belongs to $\overline{C}_1$.
We assumed that $-1$ is not pre-periodic, thus $f^{\circ m-1}(-1)$
cannot belong to the intersection $\overline{C}_0\cap\overline{C}_1$, and
the $m$-th digit in the address of $-1$ is well defined.
Denote the $m$-th digit by $\eps^*_m$.

The point $f^{\circ 2m}(-1)$ is on the boundary of $\Omega_0$.
This is the landing point of the ray $R_0(2^m\theta_0)$.
It belongs to the closure of $C_1$ or $C_0$ depending on whether
$\{2^m\theta_0\}<\theta_0$ or $\{2^m\theta_0\}>\theta_0$.
Therefore, $\eps_{2m+1}=1-\nu_m(\theta_0)$.
The point $f^{\circ 2m-1}(-1)$ is on the boundary of $\Omega_\infty$.
This is the landing point of the ray $R_\infty(2^{m-1}\theta_0)$.
It belongs to the closure of $C_0$ or $C_1$ depending on whether
$\{2^{m-1}\theta_0\}<1/2$ or $\{2^{m-1}\theta_0\}>1/2$.
Therefore, $\eps_{2m}=\theta_0[m]$.
\end{proof}

Define the following equivalence relation $\sim$ on the set of all infinite
binary sequences: $x\sim y$ if and only if one of the following formulas holds:
\begin{itemize}
\item $x=010101\dots$, $y=101010\dots$,
\item $x=w0010101\dots$, $y=w1101010\dots$,
\item $x=w0\eps^*_1\dots\eps^*_n\dots$, $y=w1\eps^*_1\dots\eps^*_n\dots$,
\end{itemize}
for some finite binary word $w$.

\begin{prop}
Let $x$ and $y$ be two infinite binary sequences.
We have $z_x=z_y$ if and only if $x\sim y$.
\end{prop}

\proof
In one direction, the proposition is obvious: if $x$ and $y$ are as described, then $z_x=z_y$.
Suppose now that $z_x=z_y$.
Interchanging $x$ and $y$ if necessary, we can write
$x=w0x'$ and $y=w1y'$ for some finite binary word $w$ (possibly empty) and infinite binary
sequences $x'$ and $y'$.
We have $z_{0x'}=z_{1y'}$.
But $z_{0x'}$ belongs to $\overline{C}_0$, whereas $z_{1y'}$ belongs to $\overline{C}_1$.
Note that the sets $\overline{C}_0$ and $\overline{C}_1$ intersect at only three points:
$\omega$, $-1$ and $a_*$.
Consider these three cases separately.

{\em Case 1.} Suppose first that $z_{0x'}=z_{1y'}=\omega$.
In this case, $x'=101010\dots$ and $y'=010101\dots$.
Indeed, if a cell lies in $C_0$ and touches the fixed point $\omega$,
then the image of this cell lies in $C_1$, and vice versa.

{\em Case 2.} Suppose that $z_{0x'}=z_{1y'}=a_*$.
In this case, it is easy to see that $x'=010101\dots$ and $y'=101010\dots$.
This follows from the fact that $f(a_*)=\omega$.

{\em Case 3}. Finally, suppose that $z_{0x'}=z_{1y'}=-1$.
Then $x'=y'=x_0$ by Proposition \ref{address_of_-1}.
$\Box$

\begin{cor}
\label{binary_model}
The Julia set of $f$ is homeomorphic to the quotient of the
space $\{0,1\}^{\N}$ of all infinite binary sequences (equipped with the product topology)
by the equivalence relation $\sim$.
Moreover, the canonical projection semi-conjugates the Bernoulli shift
with the restriction of $f$ to the Julia set.
\end{cor}

\subsection{Proof of Theorem \ref{bndry}}
Consider the two-sided lamination $2L(x_0)$, where $x_0$ is given in terms of
$\theta_0$ by the formula from Theorem \ref{bndry}.
Let us prove that the Julia set of $f$ is homeomorphic to the quotient
of the unit circle by the equivalence relation $\sim_{2L(x_0)}$, and
that the map $f$ is conjugate to the map $s_{2L(x_0)}/\sim_{2L(x_0)}$.

We can describe the equivalence relation $\sim_{2L(x_0)}$ in terms of
binary digits as follows.
Identify each point $e^{2\pi i\theta}$ on the unit circle with the binary expansion
of $\theta$, in which every second digit is replaced with its opposite.
Under this identification, the map $z\mapsto 1/z^2$ identifies with
the Bernoulli shift.

The equivalence relation $\sim_{2L(x_0)}$ is given by the following formulas:
\begin{itemize}
\item $101010\dots\sim 010101\dots$,
\item $w001010\dots\sim w11010\dots$,
\item $w0\eps^*_1\dots\eps^*_n\dots\sim w1\eps^*_1\dots\eps^*_n\dots$.
\end{itemize}
Note that the first two formulas represent identifications on the unit
circle (due to the fact that the same point on the unit circle
can correspond to different binary expansions),
and only the last formula represents the equivalence defined by
the lamination $2L(x_0)$.
The digits $\eps^*_m$ are the same as in Proposition \ref{address_of_-1}
due to Proposition \ref{binary_x_0}.

We see that the equivalence relation on binary sequences corresponding
to the relation $\sim_{2L(x_0)}$ is identical with that introduced
in Subsection \ref{ss_encoding}.
Thus both $S^1/\sim_{2L(x_0)}$ and the Julia set of the map $f$ are
identified with the quotient of the space of infinite binary sequences
by the same equivalence relation.
It follows that these two sets are homeomorphic.
Moreover, both $s_{2L(x_0)}/\sim_{2L(x_0)}$ and $f$ are represented
by the Bernoulli shift on binary sequences.
Thus the two maps are topologically conjugate.

It is easy to extend the conjugacy $(S^1/\sim_{2L(x_0)},s_{2L(x_0)})\to (J,f)$
over the gaps of the lamination $2L(x_0)$.
This finishes the proof of Theorem \ref{bndry}.

\subsection{Proof of Theorem B$^*$}
We now sketch a proof of Theorem B$^*$.
Consider a map $f\in V_2$ such that $-1\in\d\Omega_0$.
Let $\theta_0$ be the angle of the ray in $\Omega_0$ landing at the critical point $-1$.
Define the real number $y_0$ by the formula
$$
y_0=\frac 13\left(1+3\sum_{m=1}^\infty\frac{\theta_0[m]}{4^m}\right),
$$
where $\theta_0[m]$ is the $m$th binary digit of $\theta_0$.

Let $L_\Bc$ denote the basilica lamination (i.e. the lamination that models the quadratic polynomial
$z\mapsto z^2-1$).
Recall that $L_\Bc$ is a quadratic invariant lamination containing the
leaf $e^{2\pi i(1/3)}e^{2\pi i(2/3)}$ and such that this leaf has the maximal length
among all leaves in $L_\Bc$.
Consider the mating lamination $L_\Bc\cup L(-2y_0)^{-1}$.
It defines the corresponding equivalence relation on the sphere, and the quotient $S$
by this equivalence relation is also a topological sphere
(this can be deduced from a theorem of Moore \cite{Moore} that gives a
necessary and sufficient condition for a quotient of the sphere to be homeomorphic
to the sphere --- this theorem is a standard tool used to define topological matings,
see e.g. \cite{Milnor_mating}).
Let $\pi$ denote the canonical projection onto $S$.
The lamination map for $L_\Bc\cup L(-2y_0)^{-1}$ respects this equivalence relation, and,
therefore, descends to $S$.
Denote the quotient map by $\pi_*(s)$.
We would like to show that $\pi_*(s)$ is topologically conjugate to $f$.

For the map $\pi_*(s)$, we will arrange a partition into cells with exactly the
same combinatorial structure as the partition into cells of $f$.
Let $G_\infty$ be the gap of $L_\Bc$ containing the center of $G_0$ (which is a critical point for $s$),
and $G_0$ the gap containing $s(0)$ (recall that $s^{\circ 2}(0)=0$).
The open sets $\pi(G_\infty)$ and $\pi(G_0)$ are topological disks, whose boundaries intersect
at exactly one point, which is the image of the leaf $e^{2\pi i(1/3)}e^{2\pi i(2/3)}\in L_\Bc$.
Let $G_{-2}$ be the gap $-G_0$.

Let $K=\d G_0\cap S^1$.
The set $K$ is a Cantor set obtained as follows (we identify $K$ with the
corresponding subset of $\R/\Z$): from the segment $[1/3,2/3]$, we remove
two middle quarters, then do the same with the two remaining segments, etc.
Thus there is a natural parameterization of $K$ by binary sequences:
the point of $K$ corresponding to a binary sequence $\alpha_1$, $\alpha_2$, $\dots$, $\alpha_n$, $\dots$
is given by the formula
$$
\frac 13\left(1+3\sum_{m=1}^\infty\frac{\alpha_m}{4^m}\right).
$$
It is easy to see that, under this parameterization, the map $z\mapsto z^4$
(which leaves the set $K$ invariant) acts as the standard Bernoulli shift on binary sequences.

Note that the boundary of $\pi(G_0)$ is exactly $\pi(K)$.
The parameterization of $K$ by binary sequences translates into the parameterization
of $\d\pi(G_0)$ by points on the unit circle.
The parameterization of $\d\pi(G_0)$ by points on the unit circle is natural in the
sense that the map $\pi_*(s)$ restricted to $\d\pi(G_0)$ corresponds to the map
$z\mapsto z^2$ on the unit circle (this is because on binary sequences, we had the
Bernoulli shift).
The point on the unit circle parameterizing a given point $z\in\d\pi(G_0)$
will be called the {\em angle} of $z$ with respect to $\pi(G_0)$.

Consider the intersection point of $\d\pi(G_0)$ and $\d\pi(G_{-2})$.
This is a critical point for $\pi_*(s)$ (in the sense that $\pi_*(s)$ is not locally
injective near this point).
The angle of this point with respect to $\pi(G_0)$ is $\theta_0$
(this can be seen by comparing the formula for $y_0$ with the formula defining
the parameterization of $K$ by binary sequences).
It follows that the restriction of $\pi_*(s)$ to $\overline{\pi(G_\infty\cup G_0\cup G_{-2})}$
is conjugate to the restriction of $f$ to $\overline{\Omega_\infty\cup\Omega_0\cup\Omega_{-2}}$.
Now define the {\em main cell for $\pi_*(s)$} as the complement to $\overline{\pi(G_0\cup G_\infty)}$.
The {\em cells for $\pi_*(s)$} are defined as the pullbacks of the main cell under $\pi_*(s)$.
It is easy to see that the nested sequences of cells converge to points.
Thus the conjugacy between $\pi_*(s)$ and $f$ can be extended to the whole
Riemann sphere by taking pullbacks and using convergence of cells.

We see that any map $f$ on the external boundary of $M_2$ is modeled by a certain mating lamination.
We need to deduce that $f$ is a mating of the polynomial $z\mapsto z^2-1$ with some actual quadratic polynomial.
To this end, it suffices to show that the lamination $L(-2y_0)$ models a quadratic polynomial.
Note that the lamination $L(-2y_0)$ has non-renormalizable combinatorics.
Consider the parameter ray of angle $-2y_0$ in the complement to the Mandelbrot set.
For any value $c$ on the boundary of this ray, the combinatorics of the Yoccoz puzzle for
$p_c$ is the same as that for $L(-2y_0)$, in particular, $p_c$ is non-renormalizable.
It follows that $p_c$ is modeled by $L(-2y_0)$, and that $f$ is a mating of $z\mapsto z^2-1$ with $p_c$.

\subsection{Proof of Theorem B$^{**}$}
Let $f=f_a$ be such that $-1\in\d\Omega_0$.
We need to prove that $f$ is an anti-mating of $z\mapsto z^2$ with another quadratic polynomial.
We just sketch a proof skipping some details.

Consider the following family of quartic polynomials
$$
q_b(z)=bz^2(z+2)^2
$$
parameterized by a single complex parameter $b$.
For this family, $0$ is a super-attracting fixed point, and the point $-2$
is a critical point that maps to $0$.
Finally, $-1$ is a ``free'' critical point.

Clearly, if $b$ is small, then both $-2$ and $-1$ belong to the immediate basin of attraction of $0$.
Let $H$ be the hyperbolic component in the $b$-plane containing small values of $b$.
By the same methods as in \cite{Faught,Ro_puzzles}, one can show that the boundary of $H$
is locally connected.
Define the parameter ray of angle $\theta$ in $H$ as the set of all parameter values $b$
such that the critical value $b$ belongs to the interior ray of angle $\theta$ emanating from $0$.
All parameter rays in $H$ land.
Consider the landing point $b_0$ of the parameter ray of angle $\theta_0$, where $\theta_0$ is as in Theorem B$^*$.

It is not hard to see that the quartic polynomial $q_{b_0}$ is modeled by the
quartic invariant lamination $L$ from Subsection \ref{ss_L}.
From the construction of the two-sided lamination $2L(x_0)$ it is clear that
this lamination models the anti-mating of the quadratic polynomial $\sqrt{q_{b_0}}$
and $z\mapsto z^2$.

\subsection{Proof of Theorem \ref{param}}
In this subsection, we conclude the proof of Theorem \ref{param}, stating
that all external parameter rays land.
For periodic angles, this was done by Mary Rees in \cite{Rees}.
Periodic external parameter rays land at parabolic points.
For strictly pre-periodic parameter rays, the argument is essentially the same as in Proposition \ref{param_land}.
The corresponding landing points represent rational maps, for which the critical point
$-1$ is strictly pre-periodic.
Thus we can concentrate on the case of irrational angle $\theta_0$.

Consider an external parameter ray $\Rc$ of angle $\theta_0$.
Let $a$ and $a'$ be two points on the boundary of $\Rc$.
First note that, by Proposition \ref{crit_bndry} and Theorem \ref{bndry}, the maps
$f_a$ and $f_{a'}$ are topologically conjugate (since they admit the same topological model).
In particular, by Theorem B$^*$, they are matings of $z\mapsto z^2-1$ with the same quadratic polynomial.
From the Main Theorem of \cite{A-Y} it now follows that $a=a'$.

\end{document}

%% file: imsmark.tex
\def\IMSmarkvadjust{0 pt}
\def\IMSmarkhadjust{0 pt}
\def\IMSmarkhpadding{0 pt}
\def\IMSpubltext{Published in modified form:}
\def\SBIMSMark#1#2#3{
 \font\SBF=cmss10 at 10 true pt
 \font\SBI=cmssi10 at 10 true pt
 \setbox0=\hbox{\SBF \hbox to \IMSmarkhpadding{\relax}
                Stony Brook IMS Preprint \##1}
 \setbox2=\hbox to \wd0{\hfil \SBI #2}
 \setbox4=\hbox to \wd0{\hfil \SBI #3}
 \setbox6=\hbox to \wd0{\hss
             \vbox{\hsize=\wd0 \parskip=0pt \baselineskip=10 true pt
                   \copy0 \break%
                   \copy2 \break%
                   \copy4 \break}}
 \dimen0=\ht6   \advance\dimen0 by \vsize \advance\dimen0 by 8 true pt
                \advance\dimen0 by -\pagetotal
	        \advance\dimen0 by \IMSmarkvadjust
 \dimen2=\hsize \advance\dimen2 by .25 true in
	        \advance\dimen2 by \IMSmarkhadjust

%
%
  \openin2=publishd.tex
  \ifeof2\setbox0=\hbox to 0pt{}
  \else 
     \setbox0=\hbox to 3.1 true in{
                \vbox to \ht6{\hsize=3 true in \parskip=0pt  \noindent  
                {\SBI \IMSpubltext}\hfil\break
                \input publishd.tex 
                \vfill}}
  \fi
  \closein2
  \ht0=0pt \dp0=0pt
 \ht6=0pt \dp6=0pt
 \setbox8=\vbox to \dimen0{\vfill \hbox to \dimen2{\copy0 \hss \copy6}}
 \ht8=0pt \dp8=0pt \wd8=0pt
 \copy8
 \message{*** Stony Brook IMS Preprint #1, #2. #3 ***}
}

%% file: main.bbl
\begin{thebibliography}{99}

\bibitem{Ahlfors}
L. Ahlfors ``Lectures on quasiconformal mappings'',
Van Nostrand, Princeton, 1966

\bibitem{Ahmadi}
D. Ahmadi, ``Dynamics of certain rational maps of degree two'',
PhD Thesis, University of Liverpool


\bibitem{A-Y}
M. Aspenberg, M. Yampolsky, ``Mating non-renormalizable quadratic polynomials'',
preprint, {\sf http://www.arxiv.org/abs/math.DS/0610343}

\bibitem{DH}
A. Douady and J. Hubbard, ``\'Etude dynamique des polyn\^omes complex I \& II''
Publ. Math. Orsay (1984--85)

\bibitem{Faught}
D. Faught, ``Local connectivity in a family of cubic polynomials'',
PhD Thesis, Cornell University, 1992

\bibitem{Katok}
A. Katok and B. Hasselblatt, ``Introduction to the modern theory of dynamical systems''
Cambridge University Press, 1995

\bibitem{Luo}
J. Luo, ``Combinatorics and Holomorphic Dynamics: Captures, Matings and Newton's Method'',
PhD Thesis, Cornell University, 1995

\bibitem{Lyubich}
M. Lyubich, ``Six lectures on real and complex dynamics'', preprint

\bibitem{MSS}
R. Ma\~ne, P. Sud, D. Sullivan, ``On the dynamics of rational maps''. Ann. Sci. \'Ecole Norm.
Sup. (4) \textbf{16} (1983), no. 2, 193--217.

\bibitem{McMullen}
C. McMullen, ``Complex dynamics and renormalization'', Princeton University Press, 1994

\bibitem{McMullen_aut}
C. McMullen, ``Automorphisms of rational maps'', In ``Holomorphic Functions
and Moduli I'', pages 31-60. Springer-Verlag, 1988.

\bibitem{Milnor}
J. Milnor, ``Dynamics in One Complex Variable'', Third Edition, Princeton University Press, 2006

\bibitem{Milnor-QuadRat}
J. Milnor, ``Geometry and Dynamics of Quadratic Rational Maps''
Experimental Math. \textbf{2} (1993) 37--83

\bibitem{Milnor-LC}
J. Milnor, ``Local connectivity of Julia sets: expository lectures'', in
``The Mandelbrot set, Theme and Variations,''
LMS Lecture Note Series \textbf{274}, Cambr. U. Press (2000), 67--116

\bibitem{Milnor_mating}
J. Milnor, ``Pasting together Julia sets: a worked-out example of mating'',
Experimental Math \textbf{13} (2004), 55--92

\bibitem{Moore}
R.L. Moore, ``Concerning upper-semicontinuous collections of continua'',
Transactions of the AMS, \textbf{27}, Vol. 4 (1925), 416--428



\bibitem{Rees}
M. Rees, ``Components of degree two hyperbolic rational maps''
Invent. Math. \textbf{100} (1990), 357--382

\bibitem{ReesV3}
M. Rees, ``A Fundamental Domain for $V_3$'', preprint

\bibitem{Ro}
P. Roesch, ``Topologie locale des m\'ethodes de Newton cubiques'', Thesis, E.N.S, Lyon, 1997.

\bibitem{Ro_puzzles}
P. Roesch, ``Puzzles de Yoccoz pour les applications \'a allure rationnelle'', L'Enseignement Math\'ematique,
\textbf{45} (1999), p. 133--168.

\bibitem{Shub}
M. Shub, ``Endomorphisms of compact differentiable manifolds'', Amer. J. Math. \textbf{91} (1969),
175--199.

\bibitem{Sullivan}
D. Sullivan, ``Differentiable structures on fractal-like sets, determined by intrinsic scaling
functions on dual Cantor sets'', Proc. Sympos. Pure Math. \textbf{48} (1988), 15--23.

\bibitem{Sullivan1}
D. Sullivan, ``Conformal dynamical systems'', In ``Geometric Dynamics'', Springer-Verlag Lecture
Notes \textbf{1007} (1983), 725--752.

\bibitem{TanLei}
L. Tan, ``Matings of quadratic polynomials'', Erg. Th. and Dyn. Sys. \textbf{12} (1992) 589--620


\bibitem{Thurston}
W. Thurston, ``Combinatorics and dynamics of rational maps''
Princeton University and IAS Preprint. 1985

\bibitem{Tischler}
D. Tischler, ``Blaschke products and expanding maps of the circle'', Proc. Amer. Math. Soc.
\textbf{128} (1999), No 2, 621--622

\bibitem{Wittner}
B. Wittner, ``On the bifurcation loci of rational maps of degree two'',
PhD Thesis, Cornell University, 1988

\bibitem{Yampolsky}
M. Yampolsky and S. Zakeri, ``Mating Siegel quadratic polynomials'',
J. Amer. Math. Soc. \textbf{14}, No. 1, 25--78
\end{thebibliography}
